\input amstex

\loadeufm
\loadmsbm
\loadeufm

\documentstyle{amsppt}
\input amstex
\catcode `\@=11
\def\logo@{}
\catcode `\@=11
\magnification \magstep1
\NoRunningHeads
\NoBlackBoxes
\TagsOnLeft

\def \={\ = \ }
\def \+{\ +\ }
\def \-{\ - \ }

\def \b|{\big |}

\def \g1{\Gamma_1}

\def \nfp{\demo\nofrills{Proof:\usualspace\usualspace }}

\def\rarr#1#2{\smash{\mathop{\hbox to .5in{\rightarrowfill}}
 	 \limits^{\scriptstyle#1}_{\scriptstyle#2}}}

\def\larr#1#2{\smash{\mathop{\hbox to .5in{\leftarrowfill}}
	  \limits^{\scriptstyle#1}_{\scriptstyle#2}}}

\def\swarr#1#2 {\llap{$\scriptstyle #1$}  \swarrow
  	\vcenter to .5in{}\rlap{$\scriptstyle #2$}}

\topmatter
\title Automorphisms of groups and a higher rank JSJ decomposition II: The single ended case 
\endtitle
\author
\centerline{ 
Z. Sela${}^{1,2}$}
\endauthor


\footnote""{${}^1$Hebrew University, Jerusalem 91904, Israel.}
\footnote""{${}^2$Partially supported by an Israel academy of sciences fellowship.} 
\abstract\nofrills{}
The JSJ decomposition encodes the automorphisms and the virtually cyclic splittings
of a hyperbolic group. For general finitely presented groups, the JSJ decomposition
encodes only their splittings. 

In this sequence of papers we study the automorphisms 
of a hierarchically hyperbolic group (HHG) that satisfies some weak acylindricity conditions.
To study these automorphisms we construct  an object that can be viewed as a higher rank JSJ
decomposition. This higher rank decomposition encodes the dynamics of  individual automorphisms 
and the structure of the  outer automorphism group of an HHG.

\endabstract
\endtopmatter

\document

\baselineskip 12pt

The (canonical) JSJ decomposition of a  torsion-free hyperbolic group describes the structure of the (outer) automorphism group, and is
a key for understanding the dynamics of automorphisms of a hyperbolic group ([Se7],[Le]).

The construction of the JSJ was later generalized to general finitely presented groups (see [Gu-Le]). In this general setting, the JSJ encodes
all the splittings of a f.p.\ group over a given family of subgroups
(in a rather subtle way), 
but it is far from encoding the automorphism group nor
the dynamics of individual automorphisms.

In this sequence of papers we use some of the JSJ concepts, to study automorphisms of hierarchically hyperbolic groups. Hierarchically hyperbolic groups and spaces
were defined by Behrstock, Hagen and Sisto [BHS1]. The definition axiomatizes the hierarchical structure of the mapping class groups, that was defined
and studied in the work of Masur and Minsky [Ma-Mi]. Automorphisms of families of HHGs were studied earlier by  
Fioravanti [Fi], and by Casals-Ruiz, Hagen and Kazachkov ([Ca-Ka],[CHK]). 

In the first paper of the sequence, we showed how to associate a higher rank Makanin-Razborov diagram with the automorphism group of an HHG (theorem 2.8 in [Se8]).
The higher rank diagram contains a finite set of $m$-collections of (cover) resolutions, where $m$ is the number of orbits of domains of the HHG.
The higher rank diagram is constructed using a compactness argument, so it can not be canonical, but it is universal. Every automorphism of the HHG
factors through at least one of the $m$-collections in the higher rank diagram.
 
In the second paper our goal is to use the higher rank MR diagram that was constructed in the first paper and is not canonical,
 to construct objects that encode some canonical
properties and structure of both the (outer) automorphism group and the dynamics of single automorphisms of an HHG. 

To do that we start with the higher rank MR diagram of an HHG $G$, and gradually modify it to construct $m$ JSJ like decompositions
of quotients of a (characteristic) finite index subgroup $H$ of the HHG $G$. Each  decomposition is associated with an orbit of domains of the HHG $G$
under the action of the finite index subgroup $H$. We call this $m$-set of  decompositions of quotients of $H$, a $higher$ $rank$ JSJ $decomposition$ of $G$.  

Note that the strategy we are using is opposite to all the previous ones in this branch. In the case of free or hyperbolic groups, the JSJ decomposition 
is the main tool that is used to construct the Makanin-Razborov diagram. In analyzing automorphisms of HHG, 
we first constructed a Makanin-Razborov diagram using a compactness argument (as in [Ja-Se]),
which has a hierarchical structure in general. Then we gradually modify it to get the higher rank JSJ decomposition, that consists of $m$ JSJ like splittings 
with some canonical parts, and no hierarchical structure, like the JSJ decomposition of a f.p.\ group.
  
\medskip
For a one ended hyperbolic group, $\Gamma$, there exists a natural and canonical epimorphism from a finite
index subgroup of $Out(\Gamma)$
onto a direct product of the mapping class groups of the 2-orbifolds that appear in the canonical JSJ decomposition of $\Gamma$, where
the kernel is a f.g.\ virtually abelian group ([Se7],[Le]).  

In the first section of the paper we do prove a direct  generalization of these results for HHGs (theorem 1.15), 
but only in the case in which each $m$-collection in the higher rank
diagram contains at most a single non-trivial resolution. i.e., in case all the dynamics of automorphisms is concentrated or projected onto a single orbit
of domains of the HHG. 

For general HHGs $G$, in which the $m$-collections of (cover) resolutions in the higher rank MR diagram contain more than a single 
non-trivial resolutions in some of the  $m$-collections,  we did not manage to prove such a direct generalization, and
we  prove a weaker structure theorem. We replace the 
existence of a homomorphism from a finite index subgroup of $Out(G)$ into the direct product of 
finitely many mapping class groups of 2-orbifolds and outer automorphism groups
of some f.g.\ virtually abelian groups, with the construction of two groupoids that we associate with the higher rank JSJ decomposition (theorems 2.6, 3.1 and 3.2).

The two groupoids have finitely many objects. The objects in the first groupoid contain isomorphic copies of 
the finite set of QH and virtually abelian vertex groups in the $m$ decompositions in the higher rank JSJ decomposition.
The morphisms are associated with (outer) automorphisms in $Out(G)$. Each morphism is a set of  isomorphisms  between the QH 
and the virtually abelian vertex groups that are included in the domain (object) and those in the target (one isomorphism for each QH or virtually abelian 
vertex group that is associated with the object).

The objects in the second groupoid are isomorphic copies of the  edge groups in the $m$ decompositions in the higher rank JSJ decomposition.
The morphisms that are associated with (outer) automorphisms in $Out(G)$, are elements in $Z^n$, where $n$ is the total number of edges in the $m$-decompositions
in the higher rank JSJ decomposition. These elements in $Z^n$ define quasi-morphisms between $Out(G)$ and $Z^n$.

As part of the definition of a groupoid, each morphism has an inverse or a quasi-inverse, 
and it is possible to compose morphisms in a way that agrees with compositions of 
outer automorphisms in $Out(G)$. However, we were only able to associate finitely many morphisms with each
outer automorphism in $Out(G)$, and not a unique one as in the case of a hyperbolic group. 

This lack of uniqueness is a major weakness of our structure theory.  
We do not know if and where the lack of uniqueness really occurs, and we hope that in many families of HHGs it doesn't occur or it can be bypassed. e.g.,
RAAGs, (virtually) special cubulated groups etc. We see our results as an initial step, and we expect that it will be followed by quite a few
strengthenings and  refinements, as well as applications.

Although we were not able to construct a direct generalization of the structure theory for (outer) automorphisms of hyperbolic groups,
we did manage to construct a homomorphism from a finite index subgroup of $Out(G)$ into the direct sum of (finitely many) mapping class groups of modified 
2-orbifolds and general
linear groups. The 2-orbifolds are associated with the QH vertex group in the higher rank JSJ decomposition. Each 2-orbifold is constructed from the 2-orbifold
that is naturally associated with a QH vertex group by gluing 2-disks to its boundary components. Hence, the orbifold fundamental group of the closed 2-orbifold
is a (proper) quotient of the orbifold fundamental group that is naturally associated with the QH vertex group. 
The general linear groups are groups of automorphisms
of free abelian groups that are quotients of the virtually abelian vertex groups in the higher rank JSJ decomposition.

The homomorphism that we do construct gives weaker information than the epimorphism in the case of hyperbolic groups. Because of the attached disks, and
because we construct a homomorphism into the direct product and not an epimorphism. 
Nevertheless together with the groupoids that we construct, it encodes information on both
the dynamics of automorphisms and on the outer automorphism group.

We further analyze the structure of $Out(G)$ in case the higher rank JSJ decomposition 
contains no QH and no virtually abelian vertex groups. This generalizes the analysis of
the kernel of the epimorphism onto the direct sum of mapping class groups of the 2-orbifolds in the JSJ in the case of hyperbolic groups. 

In the case of hyperbolic groups this kernel is f.g.\ virtually abelian (see [Se7], [Le]). We did not manage to prove the same
result for HHGs, but we do prove that 
such $Out(G)$ are locally virtually nilpotent, by applying Gromov's theorem on groups with polynomial growth (theorem 2.8). 

\medskip 
The  construction of the groupoids from the higher rank MR diagram involves an analysis of limit quotients and sets of homomorphisms (or rather quasimorphisms)
of HHGs that encode automorphisms via the higher rank MR diagram. 
Hence, it is natural that it requires tools that were defined and constructed for analyzing varieties and more generally the first order
theory of groups, as part of the solution to Tarski's problem ([Se1]-[Se6]) and its generalizations. 
We refer the reader to some notions, objects and constructions that appear in this
sequence of papers, although we are aware that it makes the paper harder to follow. 

Throughout this sequence of papers we assume that the HHGs satisfy some weak acylindricity assumption, that appears in definition 2.1 in [Se8], and 
that the HHG is colorable. i.e., that
there exists a finite index subgroup of the HHG, under which
the orbits of the domains of the HHG are pairwise transversal. These two properties are known to hold for the mapping class groups (see [Bo], [BBF1] and [BBFS]),
and were
already assumed in the construction of the higher rank MR diagram. 

As we did in [Se8], for presentation purposes, we analyze groups that act on a product of (finitely many) hyperbolic spaces in the first two sections, 
and generalize the analysis
and the constructions
to HHGs that satisfy our weak acylindricity assumption in the third section. 

In the second paper we also assume that 
the HHGs that we analyze have a higher rank MR diagram that is single ended. This means that
the virtually abelian decompositions that are associated with the resolutions in the $m$-collections in the higher rank MR diagram contain no edges 
with finite (nor trivial) edge groups. The last assumption makes the arguments technically easier. The assumption is dropped in the next paper in the sequence
using further concepts from the work on Tarski's problem.

\vglue 1.5pc
\centerline{\bf{\S1. A higher  rank JSJ decomposition of a product I: higher rank MR}} 
\centerline{\bf{diagrams with at most a single non-trivial resolution in
each collection}}
\medskip

In section 2 of the first paper in this sequence we studied the automorphism group of a group $G$ that acts properly and cocompactly on a product
of $m$ hyperbolic spaces, and the induced action of the group on each factor is weakly acylindrical (definition 2.1 in [Se8]).
With the automorphism group of such a group we associated a higher rank Makanin-Razborov diagram (theorem 2.8 in [Se8]). 
The higher rank diagram is a finite set of $m$-collections of (cover) resolutions, where in each $m$-collection,
each cover resolution is associated with one of the factors in the product.

The higher rank MR diagram is not canonical, but it is universal. Every automorphism of the group $G$ that acts properly and cocompactly on the product,
factors through at least one of the $m$-collections of cover resolutions. 

The group $G$ that acts properly and cocompactly on the product of $m$ hyperbolic spaces, may permute the $m$ factors. 
To construct the higher rank diagram we needed to pass from the group $G$
to a finite index characteristic subgroup that we denote $H$, that fixes (setwise) the factors. The subgroup $H$ is going to be fixed for
the rest of the section.

In this section we add a simplifying assumption on the higher rank MR diagram, and obtain a higher rank JSJ decomposition.
Our goal is to prove that if each $m$-collection (of resolutions) in the higher rank MR diagram contains at most a single non-trivial
resolution (i.e., at most a single resolution with more than a single level), then the higher rank MR diagram
can be replaced by $m$ graphs of groups that are associated with the $m$ factors.
The $m$ graphs of groups have virtually cyclic edge groups, 
and their fundamental groups are quotients of the characteristic
finite index subgroup $H$. The collection of these $m$ graphs of groups encodes the dynamics of individual automorphisms in $Aut(G)$ and can
be used to study the algebraic structure of $Aut(G)$. 

In this paper we do not discuss the canonical properties of the graphs of groups that we construct and are associated with the different
factors, but we believe that the decompositions, or at least
the parts of them  that encode the
dynamics of individual automorphisms, can be made  canonical. Because of the similarity with the properties and the structure of the JSJ 
decomposition of a hyperbolic
group, we call the $m$-collection of graphs of groups that we construct, a $higher$ $rank$ JSJ $decomposition$. 

In the next section we treat products with general higher rank MR diagrams. 
In the last section we
generalize our results for products
to HHG that satisfy some further natural assumptions (that hold in the case of the 
mapping class group), that appear in section 3 in [Se8], where we construct a higher rank MR diagram for such groups. 


\medskip
Let $X$ be the product space: $X=V_1 \times V_2 \times \ldots \times V_m$, where each of the spaces $V_j$ are unbounded $\delta$-hyperbolic, for some $\delta>0$.
In the rest of this section we call the spaces $V_j$, $1 \leq j \leq m$, the factors of $X$.

Suppose that a f.g.\ group $G$ acts on $X$ isometrically, properly discontinuously and cocompactly. Let $\pi_j:X \to V_j$, $1 \leq j \leq m$, be the
natural projections. We further assume that $G$ permutes the factors $V_j$, and that 
for every index $j$, $1 \leq j \leq m$, every  $x,y \in X$ and every $g \in G$: 
$$d_{V_j}(\pi_j(x),\pi_j(y))=d_{V_{g(j)}}(\pi_{g(j)}(gx),\pi_{g(j)}(gy)).$$ Hence, a finite index subgroup of $G$ that we denote $H$, and may assume to be characteristic,
preserves the factors $V_1,\ldots,V_m$,
and acts isometrically on each of the factors $V_j$, $1 \leq j \leq m$.

Note that our assumption on the action of a group $G$ on the product space $X$ coincides with the definition of a group action on an HHS in
the special case of a product space.

The finite index subgroup $H<G$ acts isometrically on each of the factors $V_j$, $1 \leq j \leq m$. We further assume that $H$ acts weakly acylindrically on
these factors (see definition 2.1 in [Se8] for a weakly acylindrical action). 
Theorem 2.8 in [Se8] associates a finite set of collections of $m$ cover resolutions with the action of $Aut(G)$ on $H$.
In analyzing homomorphisms, the finite subset of $m$-collections of (cover) resolutions, i.e., the (higher rank) Makanin-Razborov diagram, is
the finer structure that one can construct, and in our setting in which resolutions need to be replaced by covers, it is generally not canonical (unlike
the JSJ decomposition of a hyperbolic group).
 
However, we are analyzing automorphisms, that are not just homomorphisms (or quasimorphisms). In particular, 
automorphisms are bi-Lipschitz maps,  they can be composed and they have an inverse. 
In analyzing automorphisms, our
goal is to use the finite set of  $m$-collections of  cover resolutions (theorem 2.8 in [Se8]) to construct JSJ like objects, that will encode both the dynamics of
individual automorphisms and algebraic properties of $Out(G)$.


We are going to construct the objects that encode the dynamics of automorphisms and algebraic properties of $Out(G)$ 
by gradually extending the types of cover resolutions that are contained in the 
higher rank Makanin-Razborov
diagram that is associated with $Aut(G)$. For each possible  type we use (weak) test sequences that were originally constructed to study formulas
with 2 quantifiers in [Se2], complexity of resolutions that appear in [Se4], and framed resolutions that appear in [Se6],
to construct these objects that replace the higher rank MR diagram.
We start with a rather degenerate case.

\vglue 1.5pc
\proclaim{Lemma 1.1} Suppose that $X$ is a product of $m$ unbounded (hyperbolic) factors, $V_1,\ldots,V_m$,
 $G$ is a f.g.\ group that acts properly and cocompactly on $X$, and $H$ is a characteristic finite index subgroup of $G$ that does not permute
the factors $V_1,\ldots,V_m$. 

With $Aut(G)$ and  $X$ we associate a higher rank Makanin-Razborov diagram according to theorem 2.8 in [Se8].
Suppose that an $m$-collection of cover resolutions that is part of the higher rank Makanin-Razborov diagram that is associated with $Aut(G)$, contains only
resolutions of length 1 (i.e.,
none of the resolutions contains a quotient map). Then only finitely many outer automorphisms of $G$ factor through that cover.   
\endproclaim

\nfp Let $h_1,\ldots,h_{\ell}$ be a generating set of $H$. All the $m$ cover resolutions  have length 1, so there are bounds on the lengths of the 
projections of the twisted generators:
$\varphi(h_i)$, $1 \leq i \leq \ell$, to all the spaces $V_j$, $1 \leq j \leq m$, and for all the automorphisms $\varphi \in Aut(G)$ that factor through the
collection of $m$ resolutions, after composing them with appropriate inner automorphisms.

This gives a global bound on the lengths of the elements $\varphi(h_i)$, $1 \leq i \leq \ell$, for all the automorphisms that factor through the given collection
of $m$ cover resolutions, after we compose them with inner automorphisms. Since $H$ acts properly on  $X$, 
there are only finitely many such outer automorphisms.

\line{\hss$\qed$}

\vglue 1.5pc
\proclaim{Proposition 1.2} With the assumptions of lemma 1.1, suppose that $Out(G)$ is infinite, and  
all the (finitely many) $m$-collections of cover resolutions in the higher rank MR diagram have the following properties:
\roster
\item"{(1)}" each $m$-collection contains at most one resolution that is not of a single level. i.e., at most one resolution with a quotient map.

\item"{(2)}" the resolution with a quotient map has an abelian decomposition with a single edge that has a virtually infinite abelian stabilizer, 
and no QH nor a virtually
abelian vertex group.
\endroster

Then:
\roster
\item"{(i)}"  $H$ has a  normal subgroup $N$, and a decomposition: $H=A \, *_C \, B$ or $H=A*_C$, where $N<C$, $C/N$ is virtually infinite cyclic,
and $A/N$ and $B/N$ are f.g.\

\item"{(ii)}" a finite index subgroup of $Aut(G)$ preserves $N$ and the conjugacy classes of $A,B,C$. 
Hence, it preserves the decomposition in (i).

\item"{(iii)}" $Out(G)$ is virtually cyclic. A finite index subgroup of $Out(G)$ acts as Dehn twists on the virtually cyclic decomposition of $L=H/N$:
$L=A/N \, *_{C/N} \, B/N$ or $L=A/N*_{C/N}$. 
\endroster
\endproclaim

\nfp By lemma 1.1 only finitely many (outer) automorphisms factor through an $m$-collection of  resolutions that does not contain a resolution with a quotient map.
We look at those cover resolutions from the $m$-collections of  cover resolutions that contain a quotient map. We denote them $CRes_t$, $1 \leq t \leq d$.

Let $\Lambda_t$ be the virtually abelian decompositions that are associated with the top level limit groups of the cover resolutions, $CRes_t$, $1 \leq t \leq d$. 
We denote the fundamental group of $\Lambda_t$, $L_t$.
$\Lambda_t$ contains
one or two f.g.\ vertex groups, $A_t$ and (possibly) $B_t$,  and a single edge with a virtually infinite abelian edge group, $C_t$. 
Since the action of $H$ on each of the hyperbolic spaces, $V_j$, $1 \leq j \leq m$, is assumed to be weakly acylindrical, and the edge groups are generated by
universally bounded elements,
the edge groups, $C_t$, in the virtually abelian decompositions that are associated with the top level in the resolutions $CRes_t$, have to be virtually cyclic.

With $L_t$, there is a quotient map: $\eta_t: H \to L_t$. In particular, there exist f.g.\ subgroups $\hat A_t, \hat B_t$, that together generate $H$,
and are mapped by 
$\eta_t$ to $A_t$ and $B_t$. By the construction of the higher rank Makanin-Razborov diagram,
for every automorphism $\varphi \in Aut(G)$ that factors through $CRes_t$, there is a quasi-morphism: $\mu^t_{\varphi}: L_t \to Isom(V_{j(t)})$. 

Furthermore, by the structure of the resolutions $CRes_t$, and the properness of the action of $H$ on the ambient space $X$, 
and by the existence of uniform bounds on the stretching factors of the fixed generating sets of the terminal 
limit groups of these resolutions, for any fixed $t$, there is a uniform bound on the number of  (elementwise) conjugacy classes of the images:
$\varphi(\hat A_t), \varphi(\hat B_t)$, for all the automorphisms that factor through the resolution $CRes_t$.

We look at all the possible  (boundedly many) pairs 
of (elementwise) conjugacy classes of the images of the vertex groups in the virtually cyclic  decompositions $\Lambda_t$, $1 \leq t \leq d$,
for which there is an infinite  sequence of (distinct) automorphisms
that factor through the resolutions $CRes_t$, so that the vertex groups in $\Lambda_t$ are mapped to these  pairs of conjugacy classes. 

For each infinite sequence of automorphisms $\{\varphi_s\}$, that factor through a fixed resolution, $CRes_t$, $1 \leq t \leq d$, and for which the
conjugacy classes of the images of $\varphi_s(\hat A_t)$ and $\varphi_s(\hat B_t)$ are fixed, the stable kernel of the action of $H$ on a limit tree
that is obtained from any subsequence of the sequence $\{\varphi_s\}$ is always the same. i.e., the stable kernel depends only on the conjugacy classes of the
images of the subgroups, $\hat A_t$ and $\hat B_t$, under the convergent sequence of automorphisms $\{\varphi_s\}$.
  
Since there are only finitely many resolutions $CRes_t$, and only finitely many possibilities for the pairs of conjugacy classes of the  images of the
subgroups, $\hat A_t$ and $\hat B_t$, there are only finitely many possibilities for the stable kernels of convergent sequences of automorphisms
$\{\varphi_s\}$.

Let $\{\varphi_s\}$ be a convergent sequence of automorphisms, and let $\tau \in Aut(G)$. Then $\{\varphi_s \circ {\tau}^{-1}\}$ is a convergent sequence of automorphisms
as well. Hence, every automorphism $\tau \in Aut(G)$ permutes the finite set of possible stable kernels of convergent sequences of automorphisms. Therefore,
we can pass to a finite index subgroup $U_1$ of $Aut(G)$ (or $Out(G)$), that fixes these stable kernels.

We fix a convergent sequence of automorphisms $\{\varphi_s\}$, $\varphi_s \in U_1$, that factors through a fixed resolution, $CRes_{t_0}$, where the images
$\varphi_s(\hat A_{t_0})$ are fixed, and the images $\varphi_s(\hat B_{t_0})$ are in the same conjugacy class. Let $\tau \in U_1$. Then 
for each $s$: 
$$\varphi_s(h) \ = \ \varphi_s \circ {\tau}^{-1} \circ \tau (h) \ = \ \psi_s(\tau(h))$$
where $\psi_s=\varphi_s \circ {\tau}^{-1}$.

By passing to a subsequence, we may assume that the sequence $\{\psi_s\}$ converges, factors through the same resolution, $CRes_{t'}$, and maps the
subgroups, $\hat A_{t'}$ and $\hat B_{t'}$, to fixed conjugacy classes. Since $\psi_s(\tau(\hat A_{t_0}))=\varphi_s(\hat A_{t_0})$, it follows that:
$\tau(\hat A_{t_0}) \, < \, \hat A_{t'}N$ or $\tau(\hat A_{t_0}) \, < \, \hat B_{t'}N$, where $N$ is the stable kernel. 
Similarly it follows that: $\tau(\hat B_{t_0}) \, < \, \hat B_{t'}N$
or $\tau(\hat B_{t_0}) \, < \, \hat A_{t'}N$ in correspondence.

Let $\{\varphi_s\}$ be a convergent sequence of automorphisms that factor through $CRes_{t_0}$, and let $\tau_1,\tau_2 \in U_1$ be distinct automorphisms.
Suppose that after passing to subsequences we get two convergent sequences of automorphisms,
$\{\psi_s=\varphi_s \circ {\tau_1}^{-1}\}$ and $\{\eta_s=\varphi_s \circ {\tau_2}^{-1}\}$, that factor through the same resolution $CRes_{t'}$, so that for every $s$,
$\psi_s(\hat A_{t'})$ and $\eta_s(\hat A_{t'})$ are all conjugate (elementwise), and $\psi_s(\hat B_{t'})$ and $\eta_s(\hat B_{t'})$ are also all conjugate.

Let $a_1,\ldots,a_{\ell}$ be a generating set of $\hat A_{t_0}$, and let $a'_1,\ldots,a'_{\ell'}$ be a generating set of $\hat A_{t'}$. For each $s$, 
$1 \leq i \leq \ell$, $\varphi_s(a_i)=\psi_s(\tau_1(a_i))=\eta_s(\tau_2(a_i))$. In addition, $\tau_1(\hat A_{t_0})N=\tau_2(\hat A_{t_0})N=A_{t'}N$.

Hence, $\tau_1(a_i)=w^i_1w^i_2$, where $w^i_1 \in \hat A_{t'}$ and $w^i_2 \in N$, $1 \leq i \leq \ell$. Since 
$\psi_s(\tau_1(a_i))=\eta_s(\tau_2(a_i))$ and $\psi_s(\hat A_{t'})$ is conjugate elementwise to $\eta_s(\hat A_{t'})$, it follows that there exists
some element $u \in H$ such that:
$\tau_2(a_i)=uw^i_1u^{-1} \tilde w^i_2$, $1 \leq i \leq \ell$.

$a_i={\tau_2}^{-1} \circ \tau_2(a_i)={\tau_2}^{-1}(u){\tau_2}^{-1}(w^i_1){\tau_2}^{-1}(u^{-1}){\tau_2}^{-1}(\tilde w^i_2)$, and: 
${\tau_2}^{-1} \circ \tau_1(a_i)={\tau_2}^{-1}(w^i_1){\tau_2}^{-1}(w^i_2)$. So: 
${\tau_2}^{-1} \circ \tau_1(a_i)=va_iv^{-1}w_i$, for $v \in H$, and $w_i \in N$, for $1 \leq i \leq \ell$. 

The subgroup that is generated by all the elements $\tau_2^{-1} \circ \tau_1$, for all the pairs of automorphisms $\tau_1,\tau_2 \in U_1$, for which 
for a subsequence of the automorphisms $\varphi_s$ (that factors through $CRes_{t_0}$), 
the sequences $\varphi_s \circ {\tau_1}^{-1}$ and $\varphi_s \circ {\tau_2}^{-1}$ factor through the same
resolution, $CRes_{t_0}$, and maps the subgroups $\hat A_{t_0}$ and $\hat B_{t_0}$ to conjugate subgroups, is of finite index in $U_1$, since
finitely many cosets of it cover $U_1$. 

We denote the subgroup that is generated by all these compositions ${\tau_2}^{-1} \circ \tau_1$, $U_2$. $U_2$ is of finite index in $Aut(G)$. Every
automorphism $\tau \in U_2$ preserves the stable kernel $N$, and maps $\hat A_{t_0}$ and $\hat B_{t_0}$ to their conjugates elementwise modulo the
stable kernel $N$. Hence, $U_2$ satisfies parts (i) and (ii) in the conclusion of the proposition. 

Given an automorphism $\tau \in U_2$ we can compose it with an inner automorphism, such that the composition fixes  $\hat A_0$ and maps
$\hat B_0$ to a conjugate. If $\{\varphi_s\}$ is a sequence of such automorphisms (from $U_2$), so is the sequence of compositions $\{\varphi_s \circ \tau^{-1}\}$.
Since the limit of both sequences is the Bass-Serre tree that corresponds to the decomposition: $A_{t_0}*_{C_{t_0}} \, B_{t_0}$, it follows that  each such 
automorphism $\tau$ (in $U_2$) conjugates $\hat B_{t_0}$ by an element from $\hat B_{t_0}$ mod $N$. Furthermore, since the edge group in the limit tree
is virtually cyclic, and it is its own normalizer in the limit group, modulo the kernel $N$, 
$\tau$ conjugates $\hat B_0$ by an element from the edge group that connect
the vertices $A_{t_0}$ and $B_{t_0}$. 

Therefore we get a map from the image of $U_2$ in $Out(G)$ to the edge group in the virtually cyclic decomposition of the limit group. The kernel of this map
has to be finite. Hence, $Out(G)$ has to be virtually cyclic as well, and we get part (iii).




\line{\hss$\qed$}

Proposition 1.2 proves the existence of a higher rank JSJ decomposition in case of at most a single two steps resolution, with a virtually infinite 
cyclic edge group,
and no QH nor virtually abelian vertex groups in all the $m$-collections of cover resolutions in the higher rank MR diagram.
The next theorem generalizes the construction in case there is still only at most one single two steps resolution, but the associated virtually
cyclic abelian decomposition  is only assumed not to have edges with finite (nor trivial)  stabilizers.
 
\vglue 1.5pc
\proclaim{Theorem 1.3} With the assumptions of lemma 1.1, suppose that   
in each of the  the (finitely many) $m$-collections of  cover resolutions in the higher rank Makanin-Razborov diagram, there is at most one resolution
that has two levels (i.e., a single quotient map), and all the other resolutions in each $m$-collection have a single level. 
Furthermore, assume that the virtually abelian 
decomposition that is associated with the top level of each two steps resolution contains no edges with a finite (nor trivial) edge group. 

Then:
\roster
\item"{(i)}"  $H$ admits a graph of groups decomposition that is preserved by a finite index subgroup of $Aut(G)$. If $N$ is the (pointwise) stabilizer 
of the action of $H$ on the associated Bass-Serre tree, then the associated decomposition of $H/N$ contains QH, f.g.\ virtually abelian, and rigid
vertex groups, and virtually (infinite) cyclic  edge groups. 

\item"{(ii)}" A finite index subgroup of $Out(G)$ maps into the direct sum of the mapping class groups of the QH vertex groups in the decomposition, together with
outer automorphism groups of the f.g.\ virtually abelian vertex groups (that fix their virtually abelian edge groups).
The kernel of this map is f.g.\ virtually abelian (cf. [CCV] and
[CV] for a similar theorem in the case of RAAGs, and a stronger statement in [Le] for hyperbolic groups).
\endroster
\endproclaim

\nfp If there are no two steps resolutions in all the $m$-collections, $Out(G)$ is finite by lemma 1.1. Hence, we may assume that at least in one $m$-collection
there is a two steps resolution (i.e., that $Out(G)$ is infinite). 

Given a virtually abelian abelian
decomposition with no edges with  finite edge groups $\Lambda$, we define its complexity to be a tuple of integers that are ordered lexicographically.
The tuple is composed of (in decreasing order):

\roster
\item"{(1)}" the pairs of tuples: $(-\chi(O_i),g(O_i))$ of the orbifolds that are associated with the QH vertex groups in $\Lambda$, ordered lexicographically
in non-increasing order. $\chi(O_i)$ is the Euler characteristic of the orbifold, and $g(O_i)$ is the genus of the underlying surface of each orbifold.

\item"{(2)}" the number of virtually abelian vertex groups in $\Lambda$, 
followed by a non-increasing sequence of the ranks of the maximal rank free abelian group 
in each of these
virtually abelian vertex groups. 

\item"{(3)}" the number of non-conjugate  virtually cyclic edge groups in $\Lambda$ that connect between rigid vertex groups, 
followed by a non-increasing sequence of the numbers of edges between vertex groups in $\Lambda$ that each conjugacy class stabilizes. 
\endroster
 
Let $\{CRes_t\}$, $1 \leq t \leq d$, be the two steps resolutions in the higher rank Makanin-Razborov diagram. Each automorphism in $Aut(G)$
factors through one of the resolutions in the higher rank MR diagram. For each rigid vertex group in the virtually cyclic decomposition that is associated
with $CRes_t$, we fix a f.g.\ subgroup in $H$ that is mapped onto it. Each such f.g.\ subgroup is mapped 
to a finite number of subgroups up to conjugacy, 
by all the automorphisms that factor through $CRes_t$.   

\vglue 1.5pc
\proclaim{Definition 1.4} 
We say that a sequence of automorphisms in $Aut(G)$ that factor through a resolution $CRes_t$, and converges into an action of $H$ on a real
tree, such that the virtually abelian decomposition that is associated with the action is $\Lambda_t$, is a $weak$ $test$ $sequence$ of the resolution $CRes_t$
(cf. theorem 1.18 in [Se2] for a test sequence).
\endproclaim

Since we assumed that all the cover resolutions in the higher rank MR diagram have at most two levels, the construction of the diagram guarantees
that in case of two levels  the cover resolutions have
weak test sequences that factor through them.

With each resolution $CRes_t$, we look at the finitely many tuples of possible conjugacy classes of the images of the rigid vertex groups, and the 
virtually cyclic edge groups in
$\Lambda_t$, that extend to weak test sequences of $CRes_t$. Since the limit group that is associated with a weak test sequence that factor through
$CRes_t$, is obtained from the rigid vertex groups by adding finitely many generators an relations, and since there are only finitely many possible conjugacy
classes of images of f.g.\ subgroups of $H$ that map onto each of the rigid vertex groups, there are only finitely many possibilities for the stable
kernel $N$ for all the weak test sequences that factor through a resolution $CRes_t$, $1 \leq t \leq d$.

Let $\{\varphi_s\}$ be a weak test sequence for a maximal complexity resolution $CRes_t$. For every automorphism $\tau \in Aut(G)$, a subsequence of
the sequence, $\{\varphi_s \circ {\tau}^{-1}\}$, is a weak test sequence of some maximal complexity resolution $CRes_{t'}$. Hence, $\tau$ permutes the finite
set of stable kernels of weak test sequences of maximal complexity resolutions. Therefore, a finite index subgroup $U_1 < Aut(G)$ fixes these finite
set of stable kernels elementwise. 

Since the automorphisms $\tau \in U_1$ fix the stable kernels in the limit groups that are associated with convergent sequences $\{\varphi_i \circ {\tau}^{-1}\}$,
each such automorphism restricts to an isomorphism between these limit groups. Since for any subgroup $D<H$, 
$\varphi_s(D)=\varphi_s \circ {\tau}^{-1} \circ \tau(D)$,
it follows that the restriction of an automorphism $\tau$ to an isomorphism of the limit groups of weak test sequences, permutes the conjugacy 
classes of edge groups, rigid vertex groups, virtually abelian vertex groups, and QH vertex groups. 

As in the proof of proposition 1.2, this implies that there exists a finite index subgroup $U_2<U_1$, that preserves the conjugacy classes of these edge and vertex groups
in the limit group that is associated with the weak test sequence $\{\varphi_s\}$. We denote the limit group $L$ and its virtually abelian decomposition $\Lambda$.

Hence, every automorphism $\tau \in U_2$ restricts to outer automorphisms in the mapping class groups of the QH vertex groups in $\Lambda$,
and to outer automorphisms of the virtually abelian vertex groups in $\Lambda$. Therefore, we get a map $\eta$ from
$U_2$ to the direct sum of the mapping class groups of QH vertex groups and the outer automorphism groups of virtually abelian vertex groups in $\Lambda$.

We continue with the kernel of this map $\eta$. Each automorphism from the kernel of $\eta$ preserves the stable kernel $N$ of the original sequence 
$\{\varphi_s\}$, and maps each rigid vertex group and each edge group in $\Lambda$ to a conjugate modulo the stable kernel, 
and each QH vertex group and each virtually
abelian vertex group are mapped to a conjugate modulo the kernel as well.

Therefore, the higher rank Makanin-Razborov diagram that is associated with the kernel of $\eta$ contains finitely many $m$-collections of resolutions,
where in each $m$-collection there is at most a single resolution with a (single) quotient map. Furthermore, the abelian decompositions that are associated
with these resolutions have (infinite) virtually cyclic edge groups, and no QH nor virtually abelian edge groups. 

We continue by examining a maximal complexity resolution, $CRes$, in (a maximal complexity) $m$-collection in
the higher rank MR diagram of the kernel of $\eta$, $ker(\eta)$. Let $\Lambda$ be the maximal complexity 
virtually cyclic decomposition that is associated with $CRes$.
As in the proof of part
(iii) in proposition 1.2, a finite index subgroup in $K_1< Ker(\eta)$ preserves the conjugacy classes of f.g.\ subgroups in $H$ that are mapped onto the
rigid vertex groups in $\Lambda$. Furthermore, each automorphism $\tau$ in this finite index subgroup $K_1$  acts on each vertex group in $\Lambda$ as an 
inner automorphism (i.e., it acts as an inner automorphism modulo the kernel of the action $N$).  

Therefore, as in proving part (iii) in proposition 1.2, we get a map from the image of $K_1$ in $Out(G)$ into the direct sum of the virtually cyclic edge groups in
$\Lambda$, which is a f.g.\ virtually abelian group. Hence, $Out(G)$ has to be f.g.\ and virtually abelian as well and we get part (ii) of the theorem.



\line{\hss$\qed$}


We further generalize theorem 1.3, and allow the single two step resolution in each $m$-collection to have a virtually abelian decomposition that 
includes free splittings and splittings over finite edge groups.

\vglue 1.5pc
\proclaim{Theorem 1.5} With the assumptions of lemma 1.1, suppose that   
in each of the  (finitely many) $m$-collections of  cover resolutions in the higher rank Makanin-Razborov diagram, there is at most one resolution
that has two levels (i.e., a single quotient map), and all the other resolutions in each $m$-collection have a single level. 

Then:
\roster
\item"{(i)}"  $H$ admits a graph of groups decomposition with virtually cyclic edge groups
that is preserved by a finite index subgroup of $Aut(G)$. If $N$ is the (pointwise) stabilizer 
of the action of $H$ on the associated Bass-Serre tree, then the associated decomposition of $H/N$ contains QH, f.g.\ virtually abelian, and rigid
vertex groups, and virtually cyclic (possibly finite or trivial)   edge groups. 

\item"{(ii)}" If there are QH or virtually abelian vertex groups in the graph of groups that is preserved by the finite index subgroup of $Aut(G)$,
then  a finite index subgroup of $Out(G)$ maps into the direct sum of the mapping class groups of the QH vertex groups in the decomposition, together with
outer automorphism groups of the f.g.\ virtually abelian vertex groups (that fix their virtually abelian edge groups). In this case 
the image has to be an infinite group.

\item"{(iii)}"  When the invariant virtually
cyclic decomposition of $H/N$ contains no QH and no virtually abelian vertex groups, but contains virtually infinite cyclic edge groups,
then a finite index subgroup of $Out(G)$ 
maps into a f.g.\ virtually abelian group, where the rank of this group is  the number of
virtually infinite cyclic edge groups in the invariant decomposition.
\endroster
\endproclaim

\nfp If all the two steps resolutions in the higher rank Makanin-Razborov diagram, do not have associated virtually abelian decomposition with
free decompositions nor with edges with finite edge groups, the conclusion of the theorem follows from theorem 1.3. Hence, we may assume that virtually
abelian decompositions with free splittings or with edges with finite edge groups are associated with some of the resolutions in the higher rank MR diagram.

To deal with such virtually abelian splittings we need to generalize  the definition of the complexity of a two step resolution that was previously
defined in the proof of theorem 1.3.

\vglue 1.5pc
\proclaim{Definition 1.6}    
In definition 2.1 in [Se4] we defined the $taut$ $structure$ and the $rank$ of a taut resolution over a free group. We first modify
the definition to a two step resolution in our setting.

Let $L_1 \to L_2$ be a two step resolution, and let $\Lambda$ be the virtually abelian decomposition that is associated with $L_1$.
Following definition 2.1 in [Se4] with the QH vertex groups in $\Lambda$ we can associate finitely many taut structures. With a taut structure
we can naturally associate a graph of groups, $\Lambda_T$, by cutting the orbifolds that are associated with the QH vertex groups along s.c.c.\ that
are associated with the taut structure. We assume that $\Lambda_T$ is a reduced graph of groups.
We define the Kurosh rank of the taut structure to be the number of edges with finite or trivial edge groups in $\Lambda_T$.

Given a two step cover resolution in an $m$-collection of cover resolutions from the MR diagram, we associate with it finitely many taut structures 
of its associated
virtually cyclic decomposition, by enumerating all the possible taut structures of the QH vertex groups 
in the virtually cyclic decomposition that
obtained by weak test sequences.

We continue by associating a $complexity$ of a  two steps cover resolutions that generalizes the one that appears in the proof of theorem
1.3. The complexity of a resolution is a tuple of integers, that are ordered
lexicographically. We start with the more significant terms and go down to the less significant ones.

\roster
\item"{(1)}" the maximum of the Kurosh rank over all the finitely many possible taut structures of the virtually abelian decomposition that
is associated with the resolution, that are obtained by weak test sequences.

\item"{(2)}" parts (1)-(3) in the definition of the complexity of a resolution as it appears in the proof of theorem 1.3.
\endroster
The complexity of an $m$-collection of two or one step resolutions is simply the tuple of the complexities of the $m$ resolutions in the $m$-collection 
in a non-increasing order.
\endproclaim

The set of complexities of resolutions is well-ordered. In particular, the set of the complexities of the finitely many two steps resolutions in the
higher rank MR diagram have a maximum.
Let $\{CRes_t\}$, $1 \leq t \leq d$, be the collection of two steps taut resolutions in the higher rank Makanin-Razborov diagram (that have a weak test
sequence that factors through them), that have maximal complexity.

We continue as in the proof of theorem 1.3. First, for each
rigid vertex group in the virtually cyclic decomposition that is associated
with $CRes_t$, we fix a f.g.\ subgroup of $H$ that is mapped onto it. Each such f.g.\ subgroup is mapped 
to a finite number of subgroups up to conjugacy, for all the possible weak test sequences of $CRes_t$.

With each maximal complexity taut resolution $CRes_t$, we look at the finitely many tuples of possible conjugacy classes of the images of the rigid vertex groups, 
and the 
virtually cyclic edge groups in
$\Lambda_t$, that extend to weak test sequences of $CRes_t$. Since the limit group that is associated with a weak test sequence that factor through
$CRes_t$, is obtained from the rigid vertex groups by adding finitely many generators an relations, and since there are only finitely many possible conjugacy
classes of images of f.g.\ subgroups of $H$ that map onto each of the rigid vertex groups, there are only finitely many possibilities for the stable
kernel $N$ for all the weak test sequences that factor through a maximal complexity resolution $CRes_t$, $1 \leq t \leq d$.

Let $\{\varphi_s\}$ be a weak test sequence that is taut with respect to a maximal  complexity resolution $CRes_{t_0}$. 
For every automorphism $\tau \in Aut(G)$, a subsequence of
the sequence, $\{\varphi_s \circ {\tau}^{-1}\}$, is a weak test sequence of some maximal complexity taut resolution $CRes_{t'}$.  
Hence, $\tau$ permutes the finite
set of stable kernels of weak test sequences of maximal complexity resolutions. Therefore, a finite index subgroup $U_1 < Aut(G)$ fixes this  finite
set of stable kernels elementwise. 

The rest of the proof  that there exists a  (maximal complexity) virtually cyclic splitting of a quotient  of $H$ that is invariant under a 
finite index subgroup $U_2<Out(G)$,
is identical to the argument that appears in the proof of theorem 1.3. We denote the virtually cyclic decomposition that $U_2$ preserves, $\Lambda$.

In particular, $U_2$ preserves the conjugacy classes of the QH vertex groups and the
virtually abelian vertex groups, and the edge groups that are connected to them in the invariant graph of groups that $U_2$ preserves $\Lambda$. Hence,
each automorphism in $U_2$ restricts to an outer automorphism of each of the QH vertex groups and each of the virtually abelian vertex groups
in the invariant virtually abelian decomposition. Therefore, we get a map $\eta_1$ from $U_2$ into the direct sum of the outer automorphism groups of the QH and
the virtually abelian vertex groups in $\Lambda$.

If $\Lambda$ contains no QH and no virtually abelian vertex groups, then automorphisms in $U_2$ restrict to automorphisms of the fundamental groups
of connected subgraphs of $\Lambda$, that are obtained by deleting all the edges in $\Lambda$ that have finite or trivial edge groups. By the proof of 
part (ii) in theorem 1.3, with each class of an automorphism $\tau \in U_2$ in $Out(G)$, it is possible to associate elements in each of the 
infinite virtually cyclic edge groups in $\Lambda$. This gives a map from $U_2$ into the direct sum of the edge groups with infinite virtually cyclic edge groups in
$\Lambda$, which is a f.g.\ virtually abelian group, and we get part (iii) of the theorem.

\line{\hss$\qed$}

Note that in case of free splittings or splittings over finite edge groups, the map
of $U_2$ into the direct sum of the mapping class groups of orbifolds and of the outer automorphism groups of the f.g.\ virtually abelian vertex groups,
is not guaranteed to have a f.g.\ virtually abelian kernel, because the kernel may contain automorphisms of multi-ended groups. 
It may also be the case that an MR diagram
that is associated with automorphisms that belong to the kernel of this map may have resolutions with more than two levels, and at this stage we didn't analyze such
MR diagrams.  Hence, the conclusion of
theorem 1.5 is stated in a somewhat weaker form than the conclusion of theorem 1.3.

\medskip
So far we have assumed that the $m$-collections of resolutions in the higher rank Makanin-Razborov diagram that we associated with the subgroup $H$, contain
at most a single resolution with two levels, and all the other resolutions in each $m$-collection have a single level. In the next step we still assume that 
in each $m$-collection in the higher rank MR diagram at most one (cover) resolution is not of a single level, but this one resolution can have arbitrarily many
levels. 

In this case our goal is to show that the higher rank MR diagram can be replaced by
another higher rank diagram in which in every $m$-collection of cover resolutions, all but possibly a single resolution have a single level,
and the single resolution with more than one level has at most two levels (although the resolution with two levels may not admit a weak test sequence). 
Even though it may be that the new diagram with two steps resolutions will have no weak test sequence, it will be eventually possible to 
associate with such a higher rank MR diagram a higher rank JSJ
decomposition similar to the one that was constructed in theorems 1.3 and  1.5.

Note that such a goal can not be valid when analyzing homomorphisms. In fact the multi-layer structure of a general Makanin-Razborov diagram 
that encodes homomorphisms is 
basic in the whole theory, and distinguishes the general MR diagram from a JSJ decomposition. But in analyzing automorphisms it is possible to
improve the MR diagram to have resolutions with at most two levels, and afterwards replace the diagram with a higher rank JSJ decomposition. 
This is the basic principle of our whole approach.

\vglue 1.5pc
\proclaim{Proposition 1.7} With the assumptions of lemma 1.1, suppose that   
in each of the  (finitely many) $m$-collections of  cover resolutions in the higher rank Makanin-Razborov diagram, there is at most one resolution
that has more than a single level. Suppose further that:
\roster
\item"{(i)}"  all the virtually abelian decompositions that are associated with the various levels 
of the resolutions that have more than a single level,  contain only a single edge 
with  virtually (infinite) abelian edge group, and no QH nor virtually abelian vertex groups.

\item"{(ii)}"  all the quotient maps in these resolutions  are proper quotients.

\item"{(iii)}"  a virtually abelian edge group in the virtually abelian decomposition that is associated with one of the levels
of these resolutions is not elliptic in the virtually abelian decomposition that is associated with the next level of the resolution, except
for the terminal quotient map. 
\endroster

Then:
\roster
\item"{(1)}" all the resolutions in the higher rank MR diagram that have more than a single level, have exactly two levels, and the virtually
abelian decomposition that is associated with them contains a single edge with a virtually (infinite) cyclic edge group.



\item"{(2)}"  parts (i)-(iii) of theorem 1.5 hold for $H$ and $Aut(G)$. In particular, $H$ admits a higher rank JSJ decomposition.
\endroster
\endproclaim

\nfp Let $CRes$ be one of the resolutions with more than a single level in one of the $m$-collections in the higher rank MR diagram.
Let $L$ be the limit group that is associated with the resolution $CRes_t$. Then $L$ inherits a virtually abelian decomposition from the top level
of the resolution. Each vertex group in this graph of groups decomposition inherits a virtually abelian decomposition from the virtually abelian
decomposition that is associated with the second level of the resolution, and we continue iteratively.

Since all the edge groups in the virtually abelian decompositions that are associated with the various levels of $CRes$ are virtually cyclic, all the vertex groups
and edge groups along this iterative procedure inherited decompositions are f.g. Let the terminal vertex groups 
(after the iterative decompositions) be $A_1,\ldots,A_{\ell}$.

Since in the $m$-collection that contains $CRes$, $CRes$ is the only resolution with more than a single level, for each terminal vertex group $A_i$, there are
only finitely many possible images in $H$, such that any automorphism that factors through $CRes$ sends $A_i$ to one of these images up to conjugation (elementwise).

Let $\{\varphi_s\}$ be a weak test sequence that factors through the $m$-collection that contains $CRes$, and converges to some quotient of $L$ that 
we denote $\hat L$. Let $N$ be the kernel of the quotient map from $H$ onto $\hat L$.
Since the number of images of each of the groups $A_i$ is finite up to conjugation, there are only finitely many possibilities for the kernel $N$ in $H$ for
all the weak test sequences that factor through the $m$-collection that contains $CRes$.

We define the complexity of  a  resolution that satisfies properties (i)-(iii) in the statement of the proposition, 
to be the number of quotient maps along the resolution.

We look only at resolutions with maximal complexity in the diagram, $CRes_t$, $1 \leq t \leq d$. Let $\{\varphi_s\}$ be a weak test sequence for such
a maximal complexity resolution $CRes_{t_0}$. Let $\tau \in Aut(G)$. Because of the universality of the higher rank MR diagram,
a subsequence of the sequence: $\{\varphi_s \circ {\tau}\}$ factors through one of the $m$-collections in the diagram, and in particular through
one of the 
resolutions $CRes$ in this $m$-collection.

Since the sequence $\{\varphi_s\}$ is a weak test sequence of the maximal complexity resolution
$CRes_{t_0}$, it subconverges to an action of the completion $Comp(CRes_{t_0})$ on a real tree, where the subgroups that are associated with the various levels
in $Comp(CRes_{t_0})$ act on corresponding real trees as well. The virtually abelian decompositions that are associated with these actions have  two vertices and
a single edge, where the edge group is virtually cyclic, and one of the two vertex groups is virtually abelian with a finite index subgroup
isomorphic to a free abelian group of rank 2. Furthermore, by the assumptions of
the proposition, the edge group in the abelian decomposition that is associated with a level is not elliptic in the virtually abelian decomposition that 
is associated with the next level. 

We denote the limit quotient of $Comp(CRes_{t_0})$ that corresponds to the limit action by $T$. This limit group has a tower structure similar to the
tower structure of the completion $Comp(CRes_{t_0})$, just the base group may be replaced by a proper quotient.

A subsequence of the  sequence $\{\varphi_s \circ {\tau}\}$ factors through a resolution in an $m$-collection from the higher rank MR diagram that
we denote  $CRes$. Hence, the subsequence extends to a subsequence
of quasi-actions of the completion, $Comp(CRes)$, on the corresponding hyperbolic projection space $V_j$. As in the proof of the generalized 
Merzlyakov theorem in the first section in [Se2], by possibly shorten the (extended) quasimorphisms without changing the images of $H$ under
the subsequence $\{\varphi_s \circ {\tau}\}$, the images of the quasimorphisms from $Comp(Res)$ converge into a subgroup of some closure of $T$,
$Cl(T)$ (for the notion of a closure see the first section in [Se2]).   
The closure, $Cl(T)$, is obtained from $T$ by possibly adding roots to virtually abelian vertex groups, and possibly extend the terminal limit group of
the completion. In particular, we get a map: $\eta: Comp(CRes) \to Cl(T)$.

Note that by construction the original f.g.\ group $H$ is mapped onto the limit groups that are associated with the top levels of the
completions, $Comp(CRes)$ and $Comp(CRes_{t_0})$. Since $Cl(T)$ was constructed from  subsequences of $\{\varphi_s\}$ and $\varphi_s \circ \tau\}$, 
it follows that for each $h \in H$, the image of $h$ in $Comp(CRes)$ is mapped by $\eta$ to the image of $\tau(h)$ in
$Comp(CRes_{t_0})$.

Both of the completions, $Comp(CRes)$ and $Comp(CRes_{t_0})$, are towers that are built from a base (terminal) subgroup, to which a single 
virtually abelian vertex 
group is added in each level along a virtually (infinite) cyclic edge group. Furthermore,  by our assumptions,
each of the virtually cyclic edge groups in the various
levels of the completion, $Comp(Res_{t_0})$, and
hence in the tower $T$, are not elliptic in the virtually abelian decomposition that is associated with the next level. 

$CRes_{t_0}$ is a maximal complexity resolution, i.e., it has the maximal number of levels among all the resolutions in the $m$-collections
in the higher rank MR diagram. Therefore, the existence of the map $\eta:Comp(CRes) \to Cl(T)$ that maps the base level and the top level limit groups
in $Comp(CRes)$ into the base level and the top level limit groups
in $Cl(T)$, implies that $CRes$ has to be of maximal complexity as well. 

Therefore, the automorphisms in $Aut(G)$ permute the finitely many possible kernels in $H$ of all the possible weak test sequences of maximal complexity
resolutions in the $m$-collections of the higher rank MR diagram. Hence, a finite index subgroup  $U<Aut(G)$ preserves these kernels.

Let $N$ be the kernel in $H$ that is preserved by $U<Aut(G)$. Let $C_1,\ldots,C_r$ be the virtually cyclic subgroups that are the edge groups in the
virtually abelian decompositions that are associated with the top levels of the maximal complexity resolutions in the higher rank MR diagrams. Then
the automorphisms in $U$ permute the subgroups $C_1N,\ldots,C_rN$. 

Therefore, for each fixed element $h$ in one of these subgroups, and for every $\tau \in U$, there is a global bound on the traces of the elements  $\tau(h)$
when acting on the projection space that is associated with the maximal complexity resolution. Hence, the intersection between conjugates of
the  edge groups $C_1N,\ldots,C_rN$ and $H$ have to be elliptic
in all the levels of a maximal complexity resolution. So by our assumptions, a maximal complexity resolution can not have more than two levels, and the
proposition follows.

\line{\hss$\qed$}


The next proposition generalizes proposition 1.7 to the case in which quotient maps are still proper, but edge groups in the virtually abelian decomposition that
is associated with one level need not be elliptic in the next level, and there can be more than a single edge 
in each level.  
 
\proclaim{Theorem 1.8} Suppose that only assumptions (i)-(ii) (but not necessarily assumption (iii))
in proposition 1.7 hold for the resolutions in the $m$-collections in the higher rank
MR diagram. We also allow the number of edges with infinite virtually abelian edge groups in each level to be arbitrary (and not necessarily a single edge). Then:
\roster
\item"{(1)}" The resolutions with more than a single level in the higher rank MR diagram can be replaced by a resolution with two levels,
in which each edge group is virtually (finite or infinite) cyclic, and vertex groups are either rigid or virtually f.g.\ abelian.

\item"{(2)}"  parts (i)-(iii) of theorem 1.5 hold for $H$ and $Aut(G)$. In particular, $H$ admits a higher rank JSJ decomposition.
\endroster
\endproclaim

\nfp By the construction of the higher rank MR diagram in [Se8], since we assumed that all the quotient maps in the resolutions in the $m$-collections
in the diagram are proper quotients, for each $m$-collection of resolutions in the higher rank diagram there exists a weak test test sequence of
automorphisms that factor through it. Furthermore, in case all the quotient maps are proper the resolutions are strict (definition 5.9 in [Se1]).  
 
With each resolution with more than a single level in the higher rank diagram, $CRes$, we associate its completion, $Comp(CRes)$ (see section 1 in [Se2] for
the construction of a completion). Note that by [Se2] completions can be associated with strict resolutions, and by our assumptions all the resolutions
$CRes$ are strict.

By construction, the number of levels in the completion, $Comp(CRes)$, is identical with the number of levels (or the number of quotient maps plus 1) in the
resolution $CRes$.
At this stage we gradually modify the structure of the virtually abelian decompositions that are associated with the various levels of the 
completion, $Comp(CRes)$, without changing the (limit) group that is associated with $Comp(Res)$.

If the completion has only two levels, we do not change it.
We start with the virtually abelian decomposition that is associated with the third level of the completion (counting from its bottom terminal level). 
If a virtually abelian edge group in the virtually
abelian decomposition that is associated with the third level of the resolution $CRes$ is elliptic in the virtually abelian decomposition that is associated with the
second level in $CRes$, we move the virtually abelian vertex group and the edge group that is connected to it, and both are associated with
that edge group in the third floor of the completion, to the second floor
of the completion. Note that we can move these vertex and edge groups since the edge group is assumed to be elliptic in the virtually abelian decomposition that is
associated with the second level in the resolution $CRes$.

We continue iteratively, from bottom to top. At each step $\ell$ we move all the pairs of a virtually abelian vertex group and its associated edge group 
that appear in the $\ell$-th
level of the  completion, that are associated with edge groups that are elliptic in the virtually
abelian decomposition that is associated with level $\ell-1$ in the modified completion $Comp(CRes)$ (i.e., in the completion $Comp(CRes)$ after we modified its 
bottom $\ell-1$ levels), into the $\ell-1$ level of the modified completion. And if the edge group 
is elliptic in further lower level virtually abelian decompositions, we further move the virtually abelian vertex group and the edge group that is connected
to it, to the lowest  possible level. i.e., to a level for which the edge group is not elliptic in the virtually abelian decomposition in the level below it.

At the end of the modification procedure, we did not change the group $Comp(CRes)$, but we possibly 
modified the virtually abelian decompositions of the original completion
by moving some virtually abelian vertex and edge groups to  lower levels. 

By construction, all the edge groups in the virtually abelian decompositions that are
associated with the modified completion have (infinite) virtually  cyclic edge groups. Furthermore, every edge group in a virtually abelian decomposition
which is not the lowest one, is not elliptic in the virtually abelian decomposition that is associated with the next (lower) level. 

Also, the existence of   weak test
sequences of the original completion guarantee that there exist weak test sequences for the modified completion.

\smallskip
So far we modified the completions of the resolutions so that a virtually cyclic edge group in a virtually  abelian decomposition is not elliptic in the next level
virtually abelian decomposition. However, it may still be that some of the virtually abelian vertex groups are not really essential, or that the degree
of their maximal free abelian subgroups can be further decreased. To get rid of these redundancies, and possibly reduce the number of virtually abelian 
vertex groups, or reduce their ranks, we use another procedure, that is conceptually related to the  $auxiliary$ $resolutions$ 
that were introduced in definition 2.1 in [Se5], although it is technically a different procedure.

With each of the resolutions, $CRes$,  in the higher rank MR diagram we associate its completion, $Comp(CRes)$. The construction of the completion of
a well-structured resolution appears in the first section of [Se2].
    
The limit
group $L$ (with which we started the original resolution $CRes$)  is now mapped into the completion of the resolution. We denote by $\nu(L)$ the image of
the limit group in the completion, $Comp(CRes)$ 
Also, each resolution $CRes$ has a weak test sequence, hence, each of the completions, $Comp(CRes)$ has a weak test sequence, that extends a weak test
sequence of the resolution $CRes$.

For each of the completions, $Comp(CRes)$, we define its complexity. 
The complexity of a completion, $Comp(CRes)$, is a tuple of integers ordered lexicographically in the following order:

\roster
\item"{(1)}" the number of levels in the completion.

\item"{(2)}" a tuple of integers for each virtually abelian decomposition that appears along the completion, $Comp(CRes)$, going from top to bottom.
The tuple of integers  consists of 
the number of edge groups in the virtually abelian decomposition that is associated with the  level in the completion, $Comp(CRes)$, 
followed by
the ranks of the virtually abelian vertex groups in the virtually abelian decomposition that is associated with the level in $Comp(CRes)$, ordered 
in a non-increasing order.

\item"{(3)}" a tuple of integers for each virtually abelian decomposition that appears along the completion, $Comp(CRes)$, going from top to bottom,
that consists of the number of virtually infinite cyclic edge groups in each level of the completion,
that can be conjugated into the image of the limit group $\nu(L)$. 
\endroster

Note that complexities of resolutions are well-ordered, and the higher rank MR diagram has only finitely many resolutions, so the higher rank
MR diagram contains resolutions of maximal complexity. These will be eventually used to construct the higher rank JSJ decomposition.

Let $\{\varphi_s\}$ be a weak test sequence of automorphisms in $Aut(G)$ for the resolution $CRes$. $\{\varphi_s\}$ extends to a weak test sequence of the 
completion, $Comp(CRes)$. Let $\tau \in Aut(G)$. A subsequence of $\{\varphi_s \circ \tau\}$ factors through one of the $m$-collection in the higher rank MR
diagram. Let $CRes_1$ be the non-trivial resolution in this $m$-collection (by our assumption every $m$-collection contains at most one non-trivial
resolution). 

By the construction of formal solutions and formal limit groups in [Se2], from the sequences $\{\varphi_s\}$ and $\{\varphi_s \circ \tau\}$ it is possible
to pass to a further subsequence that converges to a map:  $\nu_{\tau}:Comp(CRes_1) \to GNCl(CRes)$, where $GNCl(CRes)$ is a $generalized$ $closure$ of 
(the completion of) $CRes$. 

Recall that a closure $Cl(CRes)$ of the  completion $Comp(CRes)$ was defined in section 1 of [Se2]. It is obtained from the completion $Comp(CRes)$,
by possibly replacing virtually abelian vertex groups by finite index supergroups, and possibly changing the bottom level of the completion. A generalized
closure is obtained from a closure by possibly having elements that conjugate virtually abelian edge groups  or virtually abelian vertex groups that are not conjugate
in the completion, $Comp(CRes)$. Note that in case such conjugating elements exist the complexity of the generalized closure is strictly smaller than the 
complexity of the original completion.

Suppose that both $CRes$ and $CRes_1$ are resolutions of maximal complexity. This means that both have the same number of levels, the same number of edges in 
each level, and the same ranks of virtually abelian groups in each of the levels. Suppose further that the the generalized closure, $GNCl(CRes)$ has the same
structure as $Comp(CRes)$, i.e., no non-conjugate virtually abelian vertex or edge groups are conjugate in $GNCl(CRes)$, and that
the image of $\nu_{\tau}$ contains a conjugate of a finite
index subgroup  of each of the virtually abelian vertex group in the generalized closure $GNCl(CRes)$.

In this case we can use the same argument that was used in the proof of proposition 1.7, and deduce that the map $\nu_{\tau}$ maps edge groups
and virtually abelian vertex groups in
each level in $Comp(CRes_1)$ into edge groups and finite index subgroups of virtually abelian groups in the same level in the generalized closure, $GNCl(CRes)$.  

Hence, suppose that there exists a maximal complexity resolution in the higher rank MR diagram, $CRes$, for which there exists a weak test
sequence, $\{\varphi_s\}$, such that for every automorphism $\tau \in Aut(G)$, every subsequence of $\{\varphi_s \circ \tau\}$ that factor through a fixed 
resolution, factors through a maximal complexity resolution, $CRes_{i(\tau)}$, and the map $\nu_{\tau}:Comp(CRes_{i(\tau)} \to GNCl(CRes)$ (that is constructed
from convergent subsequences of $\{varphi_s\}$ and $\{\varphi_s \circ \tau\}$) has a graded closure, $GNCl(CRes)$, with the same structure as
$Comp(CRes)$ (i.e., no non-conjugate virtually abelian vertex or edge groups in $Comp(CRes)$ are conjugate in $GNCl(CRes)$), and
 an image that contains a conjugate of a finite index subgroup of every 
virtually abelian vertex group in the completion $Comp(Res)$.

In this case we look only at such sequences $\{\varphi_s\}$ w.r.t.\ all the maximal complexity resolutions in the higher rank MR diagram. If there exists an
edge group in a level that is above the bottom two levels in one of these maximal complexity resolutions, that has a conjugate that intersects 
the limit group $L$ in a subgroup of
finite index, then as in the proof of proposition 1.7, a finite index subgroup of $Aut(G)$ preserves the conjugacy class of the subgroup of $H$ 
that is mapped onto this edge group modulo the kernel of the action. This contradicts the assumption that the edge group is in a level that is above the bottom two levels, 
and is hyperbolic
in the level below it, which is not the bottom level.

Hence, in this case the image of $L$ does not intersect any conjugate of an edge group in a level above the bottom two levels in a subgroup of
finite index. Therefore, $L$ inherits splittings over finite edge groups from all the levels above the two bottom ones in the maximal complexity
resolutions, so the maximal complexity resolutions in the higher rank MR diagram can be replaced by resolutions with at most two levels,
and the conclusion of the theorem follows from theorem 1.5.   

\smallskip
Therefore, for the rest of the argument, we assume that for every maximal complexity resolution, $CRes$, and every weak test sequence, $\{\varphi_s\}$, that
factors through it, there is an automorphism $\tau \in Aut(G)$, such that a subsequence of $\{\varphi_s \circ \tau\}$ factors through some resolution
$CRes_1$, and after passing to convergent subsequences, either the associated map: $\nu_{\tau}:Comp(CRes_1) \to GNCl(CRes)$ 
has a generalized closure $GNCl(CRes)$ that has strictly smaller complexity than $Comp(CRes)$, or the image of $\nu_{\tau}$ 
intersects the conjugates of some virtually abelian vertex group in $GNCl(CRes)$ in subgroups of conjugates
of some fixed subgroup of infinite index of the virtually abelian vertex group in $GNCl(CRes)$. 

In this case our aim is 
to show that the maximal complexity resolution $CRes$ can be replaced by finitely many resolutions with strictly smaller complexity. Since every 
decreasing sequence of complexities of resolutions terminates after finitely many steps, an iterative modification of the higher rank MR diagram, i.e., 
an iterative replacement of maximal complexity resolutions by strictly lower complexity ones, 
concludes the proof of the theorem. 

Let $CRes$ be a maximal complexity resolution in the higher rank MR diagram, and suppose that for every weak test sequence of $CRes$, $\{\varphi_s\}$,
there exists some automorphism $\tau \in Aut(G)$, that depends on the weak test sequence, such that a subsequence of the
automorphisms, $\{\varphi_s \circ \tau\}$, factor through a resolution $CRes_1$ from the higher rank diagram, and after passing to
convergent subsequences the map: $\nu_{\tau}:Comp(CRes_1) \to GNCl(CRes)$ has either a generalized closure $GNCl(CRes)$ with strictly smaller
complexity than $Comp(CRes)$, or the image of $\nu_{\tau}$ intersects conjugates of
at least one of
the virtually abelian vertex groups along the generalized closure $GNCl(CRes)$ in subgroups of conjugates of a fixed infinite index subgroup of 
that virtually abelian vertex group.

Given a weak test sequence of $CRes$, and an automorphism $\tau \in Aut(G)$, for which the second possibility occurs.
We look at the 
(proper) image of $Comp(CRes_1)$ in $GNCl(CRes)$. We look at the highest level in $GNCl(CRes)$ that contains a virtually abelian vertex group $VA$ such that
$\nu_{\tau}(Comp(CRes_1))$ intersects conjugates of $VA$ in subgroups of conjugates of some fixed infinite index subgroup of $VA$, that we denote
$\hat VA$. 

By possibly enlarging $\hat VA$, we may assume that $\hat VA$ contains the edge group that is connected to $VA$ (even after this enlargement $\hat VA$
is of infinite index in $VA$). 
We set $\hat GNCl$ to be the subgroup of $GNCl(CRes)$, that is obtained from $GNCl(CRes)$ by replacing the virtually abelian vertex group $VA$ by its
subgroup $\hat VA$, and leaving all the other vertex groups and virtually abelian decompositions along the various levels unchanged.

By our assumptions the subsequence of the weak test sequence $\{\varphi_s\}$ that was used to construct the generalized closure, $GNCl(CRes)$, 
asymptotically factors through 
$\hat GNCl$, i.e., all the automorphisms in the subsequence factor through $\hat GNCl$ except at most finitely many. By construction, the complexity of
$\hat GNCl$ is strictly smaller than the complexity of $CRes$, since we reduced the rank of one of its virtually abelian vertex groups.

We assumed that for every weak test sequence $\{\varphi_s\}$ of $CRes$ it is possible to find an automorphism $\tau$, that enables us to construct
a generalized closure, $GNCl(CRes)$, or a modified generalized closure, $\hat GNCl$, with strictly smaller complexity. 
Clearly, there are only countable possible such generalized closures $GNCl(CRes)$, or modified generalized closures, $\hat GNCl$.
Hence, by ordering the countable set of possible generalized modified generalized closures, 
and applying the compactness argument that was used to construct the higher rank
MR diagram in section 2 of [Se8], there exists finitely many generalized and modified generalized closures: $\hat GNCl_1,\ldots,\hat GNCl_r$, 
with the following properties:

\roster
\item"{(1)}" all these generalized and  modified generalized closures have  strictly smaller complexity than
the maximal complexity of the resolutions in the original higher rank MR diagram.       

\item"{(2)}" every weak test sequence $\{\varphi_s\}$ of a maximal complexity resolution, $CRes$, in the original higher rank MR diagram, has a subsequence
that factors through one of the modified generalized closures: $\hat GNCl_1,\ldots,\hat GNCl_r$.
\endroster

In order to replace the maximal complexity resolutions in the original higher rank MR diagram by finitely many resolutions with strictly smaller
complexity, we need to further consider all the automorphisms $\sigma \in Aut(G)$, that factor only through maximal complexity
resolutions in the (original) higher rank MR diagram, and do not factor through any of the modified generalized closures, 
$\hat GNCl_1,\ldots,\hat GNCl_r$. W.l.o.g. we may consider such automorphisms that factor through a fixed maximal complexity resolution, $CRes$.

Every automorphism $\sigma \in Aut(G)$ that factors through $CRes$, extends to a homomorphism of the completion, $Comp(CRes)$, into the isometry group 
of the corresponding factor, $V_j$, of the product space $X$. We look at the set of all the extensions of the automorphisms $\sigma$ to homomorphisms of
$CRes$, for automorphisms of $\sigma$ that do do not
factor through the modified closures and through the resolutions that are not of maximal complexity in the higher rank MR diagram. 

The set of such extensions does not contain a weak test sequence of $CRes$, since otherwise a subsequence 
of the weak test sequence
factor through one of the modified closures. 

We look at sequences of extension of such automorphisms, $\sigma \in Aut(G)$, to homomorphisms of $Comp(CRes)$. We can pass to a convergent subsequence.
The convergent subsequence can not be a weak test sequence of $Comp(CRes)$, so it must converge into a proper quotient of $Comp(CRes)$, that we denote
$QC$. The map between $Comp(CRes)$ and its proper quotient $QC$ must have at least one  of the following properties:
\roster
\item"{(i)}" the rank of the image in $QC$ of some virtually abelian vertex group in $Comp(CRes)$ is strictly smaller. 

\item"{(ii)}" the images in $QC$ of two non-conjugate virtually abelian vertex groups
in $Comp(CRes)$ are conjugate in $QC$.

\item"{(iii)}"  the image in $QC$ of  an edge group in $Comp(CRes)$ can be conjugated in $QC$ into the image of the subgroup that is associated with the lower
levels in $Comp(CRes)$. In that case, in the sequence of virtually abelian decompositions that are associated with $QC$, the virtually abelian vertex
group that is connected to the image of that edge group can be pushed into a lower level.
\endroster

The limit $QC$ can be given the structure of a completion, and in all the 3 cases, the complexity of this completion has to be strictly smaller than
the complexity of $CRes$. There can be only countably many covers of such proper quotients $QC$ of $Comp(Res)$. Hence, by the compactness argument that
was used in constructing the higher rank MR diagram in section 3 in [Se8], there exist finitely many such covers, $QC_1,\ldots,QC_h$, all proper quotients of
$Comp(CRes)$, all with strictly smaller complexity than $Comp(CRes)$, with the property that all the automorphisms $\sigma \in Aut(G)$ that factor through
$CRes$, but do not
factor through the finitely many modified generalized closures, $GNCl_1,\ldots,GNCl_r$, nor through any of the resolutions in the higher rank MR diagram that are not of
maximal complexity, do factor through at least one of the quotients, $QC_1,\ldots,QC_h$.

Therefore, in case for every maximal complexity resolution, $CRes$, and every weak test sequence, $\{\varphi_s\}$, that
factors through it, there is an automorphism $\tau \in Aut(G)$, such that a subsequence of $\{\varphi_s \circ \tau\}$ factors through some resolution
$CRes_1$, and after passing to convergent subsequences, either associated map: $\nu_{\tau}:Comp(CRes_1) \to GNCl(CRes)$ is into a generalized
closure with strictly smaller complexity, or the image of $\nu_{\tau}$ 
intersects the conjugates of some virtually abelian vertex group in $GNCl(CRes)$ in subgroups of conjugates
of some fixed subgroup of infinite index of the virtually abelian vertex group in $GNCl(CRes)$, we can replace the original higher rank MR diagram
with a new higher rank diagram, such that all the automorphisms in $Aut(G)$ factor through the new diagram, and the maximal complexity
of a resolution in the new diagram is strictly smaller than the maximal complexity of a resolution in the original diagram.

The new diagram consists of the modified closures, $GNCl_1,\ldots,GNCl_r$, together with the quotient completions, $QC_1,\ldots,QC_h$, and the
resolutions in the original higher rank diagram that do not have maximal complexity. By construction, the new diagram has the universal property for
automorphisms in $Aut(G)$, and the complexities of all the (non-trivial) resolutions in it are strictly bounded by the maximal complexity of
resolutions in the original higher rank MR diagram.

Iteratively replacing the higher rank MR diagram, and strictly reducing the maximal complexity of its resolutions, we get a higher rank diagram,
for which a further complexity reduction is not possible. i.e., we get a diagram with  maximal complexity resolutions and weak test sequences
that factor through them, for which no automorphism $\tau \in Aut(G)$ enables one to apply the procedure for complexity reduction.
As we argued earlier, in this case the maximal complexity resolutions can be replaced with resolutions with at most two levels,
and the conclusion of the theorem follow by the proof of theorem 1.5.

\line{\hss$\qed$}

Theorem 1.8 proves the conclusions of proposition 1.5 in case all the $m$-collections in the higher rank MR diagram contain at most a single
non-trivial resolution, and
all the virtually abelian decompositions that are associated with the various levels in this single resolution contain only infinite virtually abelian edge groups,
and no QH nor virtually abelian vertex groups.

The next proposition generalizes theorem 1.8 by allowing the virtually abelian decompositions in the various levels of the single non-trivial resolution in
each $m$-collection to have QH and virtually abelian vertex groups.

\proclaim{Theorem 1.9} With the assumptions of lemma 1.1, suppose that   
in each of the  (finitely many) $m$-collections of  cover resolutions in the higher rank Makanin-Razborov diagram, there is at most one resolution
that has more than a single level. Suppose further that:
\roster
\item"{(i)}"  all the virtually abelian decompositions that are associated with the various levels 
of the resolutions that have more than a single level do not contain edges with finite (nor trivial) edge groups.

\item"{(ii)}"  all the quotient maps in these resolutions  are proper quotients.
\endroster

Then conclusions (1) and (2) in proposition 1.8 hold. In particular, $H$ admits a higher rank JSJ decomposition.
\endproclaim

\nfp Since all the quotient maps in the resolutions in the higher rank MR diagram are assumed to be proper quotients, each $m$-collection in the digram 
admits a weak test sequence, that restricts to a weak test sequence of the single resolution with more than a single level.

To analyze the  cover resolutions in the higher rank diagram we first need to further modify them to be $modeled$, That means to modify the completions
of 
resolutions such that some QH vertex groups along the levels of the completion,  are associated  with a subtower of the completion that is built solely on them,
and so that the modular groups that are associated with the subtower are all contained in the modular group of the QH vertex group in the bottom (base)
of the towers. Modeled resolutions
are needed to define the complexity of a resolution, while keeping the assumption that the resolution has a weak test sequence that factors through it.

\proclaim{Definition 1.10} Let $Res$ be a  well-structured   resolution 
(see  definition 1.11 in [Se2] for a well-structured resolution).  Let $Comp(Res)$ be its completion (definition 1.12 in [Se2]).   
With the completion $Comp(Res)$ we associate a $modeled$ structure.

The modeled structure of the completion is composed from a collection of subtowers (subcompletions) that are part of the completion. Each subtower is
built over a QH vertex group along the completion, which its bottom (base) level. Over the bottom QH vertex groups the subtower may contain QH and virtually
abelian vertex groups. Each QH vertex group is mapped isomorphically onto a suborbifold of the bottom QH vertex group. Each virtually abelian vertex group
is of rank 2, where the virtually cyclic 
edge group that is connected to it is mapped onto a s.c.c.\ in the bottom QH vertex group, and the virtually abelian vertex group 
is generated by its virtually cyclic edge group and a formal Dehn twist generator that is added in the construction of the completion  (see section 1 in [Se2]).

It is further required that if a QH vertex group $Q$ appears in some level of a subtower over a base QH vertex group $Q_1$, then no vertex group nor edge group 
that is not
in the subtower and appears in the completion in a level above the  bottom level of the subtower (the level of $Q_1$), is mapped (by a retraction in the tower)
to a subgroup in the subtower that intersects non-trivially  the image of the QH subgroup $Q$ by compositions of retractions in the subtower. 
Where intersects non-trivially means not in a free product of the image of conjugates of the boundary elements in $Q$.

We impose the same no intersection requirement
for the images of an edge group that is connected to a virtually abelian vertex group in the subtower.  

We denote a modeled completion, $ModComp(Res)$. It consists of the completion and the form of the subtowers that are built over QH vertex groups along it.
\endproclaim

To construct the higher rank JSJ decomposition from the higher rank MR diagram, we need to work with modeled (cover) resolutions, and not with resolutions. We
modify the cover resolutions in the higher rank diagram to be modeled, using the procedure that was used in the first part of the proof of theorem 1.8.
 
We start with a modification of the procedure from the proof of theorem 1.8. Let $CRes$ be a resolution with more than 2 levels in one of
the $m$-collections of resolutions in the higher rank MR diagram of $H$. We go along the levels of the completion, $Comp(CRes)$, and push down 
QH and virtually
abelian vertex groups, for which the virtually abelian edge groups that are connected to them  are elliptic in the lower level of the completion. We also push down 
virtually abelian vertex groups $VA$ that are connected to edges with virtually cyclic edge groups $C$, such that $VA$ are is generated by $C$ and an additional
element that is associated with Dehn twists along $C$,  and QH vertex groups $Q$, that are mapped isomorphically into a QH vertex 
group in a lower level $\hat Q$, and such that the splitting of
the group in the higher level over the virtually cyclic edge group $C$ or over the QH vertex group $Q$ is inherited from splittings of the 
group in the lower level, where this splitting of the group in the lower level
is obtained from a splitting of the QH vertex group in the lower level $\hat Q$, along a s.c.c.\ or along some suborbifold. In order to push down such a QH or a virtually
abelian vertex group, we also require that every other virtually cyclic edge group or a QH vertex group in an upper level (in the part of the
completion from the bottom up to the level that
$Q$ or $VA$ appear) is mapped to an 
elliptic subgroup in the virtually abelian decomposition of the lower level subgroup that is obtained by splitting the QH vertex group $\hat Q$ along
the isomorphic image of $Q$ or isomorphic image of $C$.  

If there are edge groups or QH vertex groups in some level of a cover resolution, that are mapped isomorphically into a QH vertex group in a lower level,
and satisfy these conditions, we push them down to be part of a tower over the QH vertex group $\hat Q$ in a lower level, which is going to be part of the modeled
structure of the obtained cover resolution.

After pushing down all the virtually abelian and QH vertex groups that satisfy the above conditions for a push down, we continue to the next upper level.
In analyzing the next level, we replace the QH vertex group $\hat Q$ that was part of some lower level of the completion,
with the QH vertex group that is obtained from $\hat Q$ by gluing QH vertex groups and virtually abelian vertex groups that are now part of the tower above it,
i.e., that is now part of its modeled structure.


We push down virtually abelian vertex groups and QH vertex groups if all of their edge groups are elliptic in a level below them or if they satisfy the conditions to
be part a tower over a QH vertex group in a lower level, taking into account that QH vertex groups in lower levels have changed according to the modeled structures
they already have.  


This variation of the procedure that was used in theorem 1.8 terminates after finitely many steps, since the number of edges and vertices in the various levels
of the completion of each cover resolution is finite, and each vertex group can move to a lower level only boundedly many times (i.e., at most the number of levels
in the completion), and can be added  to the subtower which is part of a modeled QH vertex group only once.
   
After this procedure
terminates, all the edge groups in the modified completions are virtually cyclic, and the terminal completions are all modeled. 
An edge group that is not part of a subtower over a QH vertex group in the  modeled completion, and that is not connected to a QH
vertex group in a virtually abelian decomposition that is associated with a level above
the two terminating ones in $CRes$, is not elliptic in the virtually abelian decomposition that is associated with the level below it in  of $Comp(CRes)$. Furthermore,
at least one of  the virtually cyclic 
edge groups that are connected to each of the  QH vertex groups, that appear in the bottom of each subtower over a QH vertex group in the modeled
structure of the completion,  and the QH vertex group is in a virtually abelian decomposition that is associated with a level above the two terminating ones in $CRes$,
is not elliptic in the the virtually abelian decomposition that is associated with the  level below the level in which the QH vertex group appears in $Comp(CRes)$.
 
As  in the proof of theorem 1.8, it may still be that some modeled QH or virtually abelian vertex groups in  the virtually abelian decompositions
along the modified resolutions are not really essential, 
or that the degrees
of the maximal free abelian subgroups in virtually abelian vertex groups can be further decreased. To get rid of these redundancies, 
and possibly reduce the number of modeled QH or virtually abelian 
vertex groups, or reduce their ranks, we use a procedure that generalizes the one that was used in the proof of theorem 1.8.
Conceptually, these procedures are connected to $auxiliary$ $resolutions$ that are introduced in definition 2.1 in [Se5], though they are technically different.

\smallskip
With each of the non-trivial resolutions $CRes$ in the original higher rank MR diagram, we used a modification of the procedure that was used
in the proof of theorem 1.8, to produce a modeled completed resolution, that we denote $ModComp(CRes)$. With each modeled completed resolution,
$ModComp(CRes)$, there is a map $\nu: L \to ModComp(CRes)$, that embeds the original limit group that is associated with $CRes$ into
the modeled completed resolution, and is inherited from the embedding of the limit group $L$ into the original completion, $Comp(CRes)$.
Also, each resolution $CRes$ has a weak test sequence. Hence, each of the modeled completions, $ModComp(CRes)$ has a weak test sequence, that extends a weak test
sequence of the resolution $CRes$.

As in the proof of theorem 1.8, For each of the modeled completions, $ModComp(CRes)$, we define its complexity. This generalizes the complexity of resolutions with no QH
vertex groups in the proof of theorem 1.8. To define the complexity of a modeled completion,  we regard the complexity of each modeled QH vertex group in
the modeled completion, i.e., the complexity of a QH vertex group with the subtower that is built over it, as the complexity of the QH vertex group in
the bottom of the subtower. The structure of the subtower above the QH vertex group in its bottom level, does not contribute to the complexity of
the modeled completion. 


The complexity of a modeled completion, $ModComp(CRes)$, is a tuple of integers ordered lexicographically in the following order:

\roster
\item"{(1)}" with each QH vertex group $Q$ that appears in the bottom of a modeled QH vertex group in one of the levels of the completion, $Comp(CRes)$,
we associate a pair ($-\chi(Q)$,$g(Q)$), the Euler characteristic and the genus of the associated orbifold. The highest term in the complexity
of $Comp(CRes)$ is the list of pairs that are associated with 
these  QH vertex  groups along $Comp(CRes)$, in a non-increasing lexicographical order.  

\item"{(2)}" the number of levels in the completion. In counting levels we consider only the QH vertex groups in the bottom of a subtower that is built over it, and is part
of the modeled structure of the modeled completion. We ignore the subtowers that are built over these QH vertex groups.

\item"{(3)}" a tuple of integers for each virtually abelian decomposition that appears along the completion, $Comp(CRes)$, going from top to bottom, where
the virtually abelian decompositions include only the QH vertex groups in the bottom of the subtowers that form the modeled structure of the completion.
The tuple of integers  consists of the list of pairs that are associated with the QH vertex groups (that appear in the bottom levels of modeled QH towers) 
in that level in a non-increasing lexicographical order,
followed by  
the number of edge groups in the virtually abelian decomposition that is associated with the  level in the completion, $Comp(CRes)$, 
followed by
the ranks of the virtually abelian vertex groups in the virtually abelian decomposition that is associated with the level in $Comp(CRes)$, ordered 
in a non-increasing order (edge groups and virtually abelian groups that are not part of a subtower that is built over a
modeled QH vertex group).

\endroster

As for the complexity that was used in the proof of theorem 1.8, complexities of resolutions are well-ordered. We continue by modifying the
argument that was used in proving theorem 1.8. Note that in all of our considerations modeled QH vertex groups, i.e. QH vertex groups and the subtowers
that are built over them, play the role of standard QH vertex groups. In this way we guarantee the existence of weak test sequences for the resolutions that we
construct, but introduced modeled completions, to adjust the measure of complexity, i.e., to ignore the structure of the tower that is built over a modeled QH
vertex group in the definition of the complexity, as the modular groups along the tower are all contained in the modular group of the QH vertex group in its
bottom level. 

Let $CRes$ be a non-trivial resolution in some $m$-collection in
the higher rank MR diagram. Let $\{\varphi_s\}$ be a weak test sequence of automorphisms in $Aut(G)$ for the resolution $CRes$. 
$\{\varphi_s\}$ extends to a weak test sequence of the 
completion, $Comp(CRes)$, and to a weak test sequence of the modeled completion, $ModComp(CRes)$. 
Let $\tau \in Aut(G)$. A subsequence of $\{\varphi_s \circ \tau\}$ factors through one of the $m$-collections in the higher rank MR
diagram. Let $CRes_1$ be the non-trivial resolution in this $m$-collection (by our assumption every $m$-collection contains at most one non-trivial
resolution), and let $ModComp(CRes_1)$ be its modeled completion.

The modeled completion has a weak test sequence that factors through it, and it was used to define the complexity of a (completed) resolution. However, to
compare between different completed resolutions that a sequence of values factors through, we need to further reduce the structure of a modeled completion,
to a $reduced$ ($modeled$) $completion$.

\proclaim{Definition 1.11} Let $Res$ be a  resolution, and let $ModComp(Res)$ be its modeled completion. The $reduced$ $modeled$ $completion$,
that we denote $RedComp(Res)$  is the completion that is obtained from the modeled completion, $ModComp(Res)$, by collapsing each subtower that is built
over a QH vertex group and is part of the modeled structure in $ModComp(Res)$ to the QH vertex group in the bottom of the subtower.

Note that any value that factors through the model completion, factors through the reduced modeled completion. On the other hand, a weak test sequence of the 
modeled completion is not a weak test sequence of the reduced modeled completion.   

\endproclaim

Note that the complexity of a modeled resolution can be defined using the complexity of its reduced resolution. To analyze and modify the higher rank MR diagram
we will also need $framed$ modeled resolutions. Framed resolutions were originally defined in definition 5 in [Se6], and were used in the  procedure for
quantifier elimination over a free group.

\proclaim{Definition 1.12} Let $Res$ be a  resolution, and let $ModComp(Res)$ be its modeled completion. A framed resolution of the modeled completion
$ModComp(Res)$, that we denote $FrmComp(Res)$, is a completion that is obtained from $ModComp(Res)$ by replacing some of the QH vertex groups along
$ModComp(Res)$ with a finite index supergroups, adding finite order roots to some of the virtually abelian vertex groups, and possibly adding elements to the 
bottom level of the modeled completion $ModComp(Res)$.

Since a framed resolution is a completion, the retractions of the modeled completion, $ModComp(Res)$, need to be modified to include the
(images of the)  finite index supergroups of the original QH vertex groups. Also if a QH vertex group $Q$  in a subtower that is part of the modeled structure
of $ModComp(Res)$ is replaced by a finite index supergroup, then all the QH vertex groups that are below $Q$ in that subtower are replaced by finite index
supergroups. In particular, in that case the QH vertex group in the bottom of the subtower is replaced by a finite index supergroup.

If some of the QH vertex groups in a subtower are replaced by finite index supergroups (in that case that includes
the QH vertex group in the bottom of the subtower) and some remains unchanged, the modeled structure is changed, since the modular groups of some QH vertex groups
are not contained (after the modification) in the modular group of the QH vertex group in the bottom of the original subtower.

Hence, a framed resolution of a modeled resolution may have a different (more refined) modeled structure than the original modeled resolution. However, in case
a QH vertex group is replaced by a proper finite index supergroup, the complexity of the framed resolution is strictly smaller than the complexity
of the original modeled resolution.   
\endproclaim

As in the proof of theorem 1.8, by the construction of formal solutions and formal limit groups in [Se2], 
from the sequences $\{\varphi_s\}$ and $\{\varphi_s \circ \tau\}$ it is possible
to pass to a further subsequence that converges to a map:  $\nu_{\tau}:RedComp(CRes_1) \to FCl(CRes)$. $FCl(CRes)$ is a framed generalized closure 
that is obtained from a closure of the modeled
completion $ModComp(CRes)$, $Cl(CRes)$, by possibly replacing some QH vertex groups and virtually abelian vertex groups in the closure 
with a finite index supergroups that contain them, and possibly having elements that conjugate virtually abelian or QH vertex groups that are not conjugate in
the modeled completion, $ModComp(CRes)$.

\proclaim{Lemma 1.13} Suppose that $CRes$ is a resolution of maximal complexity in the higher rank MR diagram and $FCl(CRes)$ is a framed (generalized) closure
of it. Then 
either:
\roster

\item"{(1)}" two non-conjugate virtually abelian or QH vertex groups or two non-conjugate edge groups in $ModComp(CRes)$ are mapped to conjugate ones
in $FCl(CRes)$.

\item"{(2)}" Some QH vertex groups in the framed closure $FCl(CRes)$ properly contain QH vertex groups in the modeled resolution $ModComp(CRes)$.
In that case, some QH vertex groups that appear in the bottom of (modeled) subtowers in $ModComp(CRes)$ are contained in proper finite index supergroups
in the framed closure $FCl(CRes)$. The framed closure $FCl$ inherits a modeled structure from that of the modeled (completed) resolution,
$ModComp(CRes)$, and the modeled structure of $FCl$ may be further refined (see definition 1.12). This means that modeled subtowers in $ModComp(CRes)$ may be
divided into finitely many modeled subtowers in (the modeled completion structure of) $FCl(CRes)$. In particular, since at least one QH vertex group in the bottom
of a subtower in $ModComp(CRes)$ is replaced by a proper finite index supergroup in $FCl(CRes)$, the complexity of the framed closure $FCl(CRes)$ 
is strictly smaller than the complexity of
the (maximal complexity) modeled resolution $ModComp(CRes)$.
\endroster

Suppose that the framed closure, $FCl(CRes)$, does not satisfy part (1) nor (2). In that case, after possibly modifying the structure of $FCl(CRes)$ and
the map $\nu_{\tau}$,  either:
\roster
\item"{(3)}" $\nu_{\tau}(RedComp(CRes_1))$ intersects some QH vertex group and all the subtower over it in the modeled closure  $FCl(CRes)$ trivially. 
i.e. in products of conjugates of some of the boundary components in the bottom QH vertex group.

\item"{(4)}" $\nu_{\tau}$ maps  the QH vertex groups in $RedComp(CRes_1)$ isomorphically into modeled  QH vertex groups in $FCl(CRes)$, and every modeled
QH vertex group in
$FCl(CRes)$  has a QH vertex group in $RedComp(CRes_1)$ that is mapped non-trivially into (a conjugate of) 
the subtower that is built over it by the map $\nu_{\tau}$.
There exists some virtually abelian vertex group $VA$ in $FCl(CRes)$, such that the
intersection between $\nu_{\tau}(Comp(CRes_1))$ and conjugates of $VA$ are contained in conjugates of some subgroup of infinite index: $\hat VA<VA$.    

\item"{(5)}" $CRes_1$ is a resolution of maximal complexity in the higher rank MR diagram, $\nu_{\tau}$ maps the QH vertex groups in $Comp(CRes_1)$
isomorphically into conjugates of subtowers that are built over distinct QH vertex groups in $FCl(CRes)$, and virtually abelian vertex groups in $Comp(CRes_1)$ 
into conjugates of 
finite index subgroups of distinct virtually abelian vertex groups in $FCl(CRes)$. 
\endroster
\endproclaim

\nfp 
Let $\nu_{\tau}:RedComp(CRes_1) \to FCl(CRes)$. The framed closure $FCl(CRes)$ is obtained from a closure of the modeled completion, 
$ModComp(CRes)$, by either conjugating two virtually abelian or QH vertex groups or virtually abelian edge groups (part (1)), or by 
possibly replacing some QH vertex groups with finite index supergroups that contain them and changing the retraction maps between the
levels of the framed closure to fit with this replacement (part (2)). If there is such a QH vertex group $Q$ in
$ModComp(CRes)$, then there is such a QH vertex group in the bottom of a subtower that is associated with the modeled structure of $ModComp(CRes)$,  
 and this QH vertex group is replaced  in
$FCl(CRes)$ by a  QH vertex group that properly contains it. 
In this case part (2) in the theorem holds. In both cases (1) and (2) the complexity of the new modeled structure of $FCl(CRes)$ 
is strictly smaller than the complexity of
$ModComp(CRes)$.

In the sequel suppose that parts (1) and (2)  in the theorem do not hold.
If the maximal complexity resolution, $CRes$, contains no QH vertex groups, the lemma follows from 
theorem 1.8. If $CRes_1$ contains no QH vertex groups, the lemma follows from the proof of theorem 1.8.

Suppose that $CRes_1$, and hence the maximal complexity resolution $CRes$, contains QH vertex groups. We start by collapsing the subtowers in the
modeled structure of the framed closure $FCl(CRes)$ 
(if parts(1) and (2) do not occur, $FCl(CRes)$ inherits the modeled structure of $ModComp(CRes)$), to get a $reduced$ framed closure
$RedFCl(CRes)$.

We further compose  $\nu_{\tau}$ with the retractions
in the modeled subtowers in the framed closure $FCl(CRes)$  to get a map:
$\hat \nu_{\tau}:RedComp(CRes_1) \to RedFCl(CRes)$. 

For presentation purposes, suppose
first that $CRes_1$ contains a QH vertex group $Q_1$ in its terminal virtually abelian decomposition. i.e.,  all the boundary components of $Q_1$
are contained in the terminal bounded group in $CRes_1$. Note that $Q_1$ must be in the bottom of a subtower that is part of the modeled structure of 
$ModComp(CRes_1)$, hence, it is part of $RedComp(CRes_1)$.

By construction, $\hat \nu_{\tau}$ maps the boundary elements of $Q_1$ into the bounded subgroup in $RedFCl(CRes)$. Suppose that $\hat \nu_{\tau}(Q_1)$ intersects
conjugates of some QH vertex groups or some virtually abelian vertex groups  along the levels of $RedFCl(CRes)$ non-trivially. i.e., 
in the case of QH vertex group, in subgroups that are not conjugate to free products of conjugates
of boundary elements. Suppose that  $Q_u$ is a QH vertex group in $RedFCl(CRes)$, that intersects non-trivially a conjugate of $\hat \nu_{\tau}(Q_1)$,
and $Q_u$ is in the  highest level in $RedFCl(CRes)$ for which there are such $QH$ or virtually abelian vertex groups.

Since the virtually abelian decompositions 
along the resolution
$CRes_1$ contains no edges with finite edge groups ($CRes_1$ is of minimal rank following the definition in section1 in [Se4]), 
a suborbifold of $Q_1$ is mapped by $\hat \nu_{\tau}$ onto a finite index subgroup of a conjugate of $Q_u$. 
Therefore,
the complexity (pair) of $Q_u$ is bounded by the complexity of $Q_1$, with equality if and only if  $Q_u$ is a QH vertex group in the terminal virtually abelian
decomposition of $RedFCl(CRes)$, and $\hat \nu_{\tau}$ maps $Q_1$ isomorphically onto a conjugate of $Q_u$. 
 
Let $i$ be the highest level in the reduced framed closure, $RedFCl(CRes)$, for which $\nu_{\tau}(Q_1)$ intersects non-trivially a conjugate of a 
QH or a virtually abelian vertex group in that closure. Let $RedFCl(CRes)_{i+1}$ be the subgroup which is the part of the closure, $RedFCl(CRes)$,
that contains all the levels up to level $i+1$ (the level below level $i$).  Let $\eta_{i+1}$ be the retraction: $\eta_{i+1}:RedFCl(CRes) \to RedFCl(CRes)_{i+1}$.

We look at the image of $Q_1$ under the composition: $\eta_{i+1} \circ \hat \nu_{\tau}$. The conjugacy classes of the
images of boundary components of $Q_1$ do not change by the composition with
the retraction $\eta_{i+1}$. Suppose that the image of $Q_1$ under the composition intersects non-trivially 
a conjugate of a QH vertex group $Q^{i+1}$ in level $i+1$ 
of $RedFCl(CRe)$, which is
the highest level in the part of $RedFCl(CRes)$ that is the image of the retraction $\eta_{i+1}$. 

By the same arguments that applied to the QH vertex group
$Q_u$, the complexity of $Q^{i+1}$ is bounded by the complexity of $Q_1$, with equality 
if and only if  $Q^{i+1}$ is a QH vertex group in the terminal virtually abelian
decomposition of $RedFCl(CRes)$, and $\eta_{i+1} \circ \hat \nu_{\tau}$ maps $Q_1$ isomorphically onto a conjugate of $Q^{i+1}$. 

We continue iteratively, by composing with retractions to lower levels. By the same arguments we conclude that any QH vertex group that has a conjugate
that is intersected non-trivially
by the image of $Q_1$ under the composition of a retraction and the map $\hat \nu_{\tau}$, has complexity that is bounded by the complexity of $Q_1$, with equality
if and only if the QH vertex group is in the bottom level of $RedFCl(CRes)$, and the composition of the retraction and $\hat \nu_{\tau}$ maps $Q_1$ isomorphically
onto that QH vertex group.

\smallskip
The bounds on the complexities of the QH vertex groups in $RedFCl(CRes)$ that intersect non-trivially conjugates of the images of a
QH vertex group $Q_1$ in the bottom virtually abelian decomposition of $RedComp(CRes_1)$, under compositions of retractions of $RedFCl(CRes)$ with the map
$\hat \nu_{\tau}$, enable us to analyze the image of $\hat \nu_{\tau}$.

If the complexities of all the QH vertex groups in $RedFCl(CRes)$ that have conjugates that are intersected non-trivially by $\hat \nu_{\tau}(Q_1)$ and its 
images under the retractions, 
have strictly smaller complexity than $Q_1$,
we gain in the ambient 
complexity. i.e., $Q_1$ covers only QH vertex groups in $RedFCl(CRes)$ that have strictly smaller complexity.

If there is a QH vertex group $Q$ in $RedFCl(CRes)$ that has a conjugate that is intersected non-trivially by $\hat \nu_{\tau}(Q_1)$ or its images under retractions,
and $Q$
has the same complexity as $Q_1$, $Q$  has to be in the terminal virtually abelian decomposition of $RedFCl(CRe)$, and the composition of the map $\hat \nu_{\tau}$
with (possibly) retractions of $RedFCl(CRes)$ have to map $Q_1$  onto a conjugate 
of $Q$ in $RedFCl(CRes)$. 

If there are two QH vertex groups $Q_1$ and $Q_2$
in the terminal virtually abelian decomposition of $CRes_1$, such that the images of them under a composition of a retraction with $\hat \nu_{\tau}$ intersect conjugates
of the same QH vertex group $Q$ in the terminal virtually abelian decomposition of $FCl(CRes)$, we also gain in the ambient complexity.

\proclaim{Definition 1.14 (cf. definition 1.6 in [Se4])} A QH vertex group $Q$ in the terminal virtually abelian decomposition of $RedFCl(CRes)$ 
is called $surviving$ $orbifold$ if:

\roster
\item"{(i)}"  there exists 
a QH vertex group $Q_1$ in the terminal virtually abelian decomposition of $RedComp(CRes_1)$ such that
the composition of a retraction and $\hat \nu_{\tau}$ maps $Q_1$ isomorphically onto a conjugate of $Q$.

\item"{(ii)}" the images of all the other QH vertex groups in the terminal level of $RedFCl(CRes_1)$ under the composition of a retraction and $\hat \nu_{\tau}$,
intersect every conjugate of $Q$ trivially (i.e., in a free product of conjugates of boundary subgroups of $Q$).
\endroster
\endproclaim

Suppose that $Q$ is a surviving orbifold, and let $Q_1$ be the unique QH vertex group in the terminal virtually abelian decomposition
of $CRes_1$, that is mapped by a composition of (possibly) a retraction and $\hat \nu_{\tau}$ isomorphically onto a conjugate of $Q$. 
We continue in a similar way to what was done
in section 1 in [Se4] for bounding the complexity of quotient resolutions.

We  rearrange QH and virtually abelian 
subgroups along the levels of the reduced framed closure, $RedFCl(CRes)$, that are mapped by retractions isomorphically onto suborbifolds of $Q$,
and modify accordingly the map $\nu_{\tau}$. We do that to guarantee that for each surviving orbifold $Q$
in the terminal virtually abelian decomposition of $RedFCl(CRes)$,
the unique QH vertex group $Q_1$ that was mapped isomorphically onto a conjugate of $Q$  by a composition of a retraction and the map $\hat \nu_{\tau}$,
will be mapped isomorphically onto $Q$ by the modification of $\hat \nu_{\tau}$ (with no composition with a retraction). This will imply that a QH vertex group
$Q_1$ in $RedComp(CRes_1$ that is mapped onto a conjugate of a surviving orbifold $Q$, intersects trivially conjugates of all the other QH vertex groups
in $RedFCl(CRes)$. Hence, the contribution of $Q_1$ to the ambient complexity of $CRes_1$ will be the same as the contribution of the surviving orbifold $Q$ to
the complexity of $RedFCl(CRes)$,  that is identical to the complexity of $ModComp(CRes)$.

Let $Q^i,Q^{i+1},\ldots,Q^{t-1}=Q$, be the sequence of (not necessarily connected, some of which possibly empty) 2-orbifolds and virtually
abelian vertex groups in $Cl(CRes)$, from level $i$ down to the QH vertex group $Q$ in the terminal virtually abelian decomposition in level $t-1$, 
such that the image of $Q_1$ under a composition of the corresponding
retraction with $\hat \nu_{\tau}$ intersects non-trivially a conjugate of them. 

Since $Q_1$ is mapped isomorphically onto a conjugate of the surviving orbifold $Q$ 
by a composition of a retraction and $\hat \nu_{\tau}$, suborbifolds of $Q_1$ are mapped isomorphically onto conjugates of the QH vertex groups in the sequence:
$Q^i,\ldots,Q^{t-1}$, and s.c.c.\ in $Q_1$ are mapped onto finite index subgroups of conjugates of edge groups in virtually abelian
vertex groups along this sequence. Furthermore, the final retraction maps QH vertex groups in $Q^i,\ldots,Q^{t-2}$ isomorphically onto conjugates of suborbifolds of
$Q$, and finite index subgroups of edge groups of virtually abelian vertex groups onto s.c.c.\ in conjugates of $Q$. 
Since $Q$ is a surviving orbifold, the images of all the other QH vertex groups in the terminal level of $ModComp(CRes_1)$, under compositions of
appropriate retractions with $\hat \nu_{\tau}$, intersect trivially conjugates of all the QH and virtually abelian vertex groups in the sequence:
$Q^i,\ldots,Q^{t-1}$.

These observations allow us to modify the structure of the reduced framed closure $RedFCl(CRes)$ and the map $\hat \nu_{\tau}$ in a similar way to the 
modification of quotient resolutions
(in the minimal rank case)
in the first sections of [Se4] and [Se5].

We modify the reduced framed closure $RedFCl(CRes)$ by inverting the order of the levels of the QH and virtually abelian vertex groups that appear in the sequence:
$Q^i,\ldots,Q^{t-2}$, leaving $Q^{t-1}=Q$ unchanged. i.e., we set it to be: $Q^{t-2},\ldots,Q^i,Q^{t-1}=Q$. 
This order inversion is completed by changing accordingly the map $i\hat \nu_{\tau}$,
that after modification maps $Q_1$ isomorphically onto a conjugate of $Q$ (without composing it with a retraction). 
For the detailed modification of the map $i\hat \nu_{\tau}$ see section 1 in [Se4].
We denote the modified map $\tilde \nu_{\tau}$, and the modified reduced framed closure, $\tilde RedFCl(CRes)$. 

Note that the images of the other QH and virtually abelian vertex groups in the terminal virtually abelian decomposition  in  $RedComp(CRes_1)$ are not effected by
the modification of $\hat \nu_{\tau}$, as these images intersect the QH and virtually abelian vertex groups in the sequence $Q^i,\ldots,Q^{t-1}=Q$ trivially.
The images of QH and virtually abelian vertex groups in higher level virtually abelian decompositions in $RedComp(CRes_1)$ may be effected by the modification
of the map $\hat \nu_{\tau}$, and the same holds for the image of $L$, $\nu(L)$, in the modified closure $\tilde RedFCl(CRes)$ (see section 1 in [Se4] for 
these modifications).

We perform these modifications of the reduced framed closure $RedFCl(CRes)$, and of the map $\hat \nu_{\tau}$, 
for all the surviving orbifolds. Note that because of trivial intersections,
these modifications can be conducted in parallel. Also note that because the edge groups that are connected to the virtually abelian vertex groups in the terminal 
virtually abelian decomposition of
$RedComp(CRes_1)$ do all belong to the terminal bounded subgroup in $RedComp(CRes_1)$, $\hat \nu_{\tau}$ and the modified map $\tilde \nu_{\tau}$ map these virtually
abelian vertex groups into  the terminal bounded subgroup in $\tilde RedFCl(CRes)$, or into conjugates of virtually abelian vertex groups in the terminal
virtually abelian decomposition in $\tilde RedFCl(CRes)$.

\smallskip
We continue iteratively by climbing along the levels of $RedComp(CRes_1)$, 
and possibly further modifying the structure of the reduced framed closure 
$\tilde RedFCl(CRes)$ and the map $\tilde \nu_{\tau}$.
After completing the modification at level $i$, $2 \leq i \leq t-1$, we set $M_i$ to be the subtower of the  reduced framed closure, $\tilde RedFCl(CRes)$, 
that contains the QH and 
virtually abelian vertex groups in $RedFCl(CRes)$ that intersect non-trivially images (under $\tilde \nu_{\tau}$) 
of QH and virtually abelian vertex groups in the
levels $i$ and below in the reduced completion $RedComp(CRes_1)$. 
By construction,  $M_i$  is a subtower that is closed under the retractions of
the ambient modified closure $\tilde RedFCl(CRes)$.

Note that the images under $\tilde \nu_{\tau}$ of all the edge groups in the virtually abelian decomposition in level $i-1$ in $RedComp(CRes_1)$, including all the 
boundary elements of QH vertex groups in this virtually abelian decomposition, are contained in $M_i$. Hence, if a QH vertex group $Q_{i-1}$
in the virtually abelian decomposition that is associated with level $i-1$ in $RedComp(CRes_1)$, and $Q$ is a QH vertex group in the highest possible level in 
$RedFCl(Res)$
that
is not in $M_i$ and intersects non-trivially a conjugate of $\tilde \nu_{\tau}(Q_{i-1})$, then the complexity of $Q$ is bounded by the
complexity of $Q_{i-1}$.  Furthermore, equality in the complexities of $Q$ and $Q_{i-1}$ occurs if and only if $Q$ is a QH vertex group in some level of 
$RedFCl(CRes)$,
and all the QH and virtually abelian vertex groups in the levels below it, that contain conjugates of its boundary components, are in the subtower of
$RedFCl(CRes)$ that is associated with $M_i$. 

We look at compositions of $\tilde \nu_{\tau}$ and retractions of the closure, $RedFCl(CRes)$, precisely as we did in analyzing QH vertex groups
in the terminal
virtually abelian decomposition of $RedComp(CRes_1)$. 

The analysis is complete when we analyze the highest level in $RedComp(CRes_1)$. Suppose that part (3) in the conclusion of the theorem
does not hold. Since $CRes$ was assumed to be of maximal complexity, the definition of the
complexity of a modeled completion implies that if all the QH vertex groups in $RedFCl(CRes)$ 
intersect non-trivially conjugates of the image of
$RedComp(CRes_1)$ under $\tilde \nu_{\tau}$,   then every $QH$ vertex group in $RedComp(CRes_1)$
is mapped isomorphically by $\tilde \nu_{\tau}$  onto a conjugate of a QH vertex group in 
$RedComp(CRes)$, and these QH vertex groups are distinct. 

In that case we analyze the image under $\tilde \nu_{\tau}$ of the virtually abelian vertex groups in 
$RedComp(CRes_1)$, precisely as we did in the proof of theorem 1.8. If there exists a virtually abelian vertex group in one of the levels of $RedFCl(CRes)$
that intersects $\tilde \nu_{\tau}(RedComp(CRes_1))$ in a subgroups that can be conjugated into some fixed infinite index subgroup of the virtually abelian 
vertex group in $RedFCl(CRes)$ we get part (4) of the lemma. If there isn't such a virtually abelian vertex group, then by the proof of theorem 1.8,
$\tilde \nu_{\tau}$ maps the virtually abelian vertex groups in $RedComp(CRes_1)$ isomorphically onto finite index subgroups of distinct virtually abelian
vertex groups in $RedFCl(CRes)$. In particular, $CRes_1$ has to be a maximal complexity resolution, and we get part (5) of the lemma.  
   
\line{\hss$\qed$}

We continue the proof of theorem 1.9 along the lines of the argument that was used to prove theorem 1.8.
Let $CRes$ be a maximal complexity resolution in the higher rank MR diagram, and suppose that there exists a weak test sequence,
$\{\varphi_s\}$ in $Aut(G)$, such that for every automorphism $\tau \in Aut(G)$, all the convergent pairs of subsequences,
$\{\varphi_s\}$ and $\{\varphi_s \circ \tau\}$, converge into maps: $\nu_{\tau}:RedComp(CRes_1) \to RedFCl(CRes)$, that satisfy part (5) in lemma 1.13.

In this case, we look at the set of all the weak test sequences that satisfy this condition w.r.t.\ all the maximal complexity resolutions
in the higher rank MR diagram. We argue as in the proofs of proposition 1.7 and theorem 1.8.  
If there exists an
edge group in a level that is above the bottom two levels in one of the reduced model completions of these maximal complexity resolutions, $RedComp(CRes)$,
and the edge group
  has a conjugate that intersects 
the image of the limit group $L$ in $RedComp(CRes)$  in a subgroup of
finite index, then  a finite index subgroup of $Aut(G)$ preserves the conjugacy class of the subgroup of $H$ 
that is mapped onto this edge group. This contradicts the assumption that the edge group is in a level that is above the bottom two levels, and is hyperbolic
in the level below it, which is not the bottom level.

Hence, in this case the image of $L$ does not intersect any conjugate of an edge group in a level above the bottom two levels in a subgroup of
finite index. Therefore, $L$ inherits splittings over finite edge groups from all the levels above the two bottom ones in the maximal complexity
resolutions, so the maximal complexity resolutions in the higher rank MR diagram can be replaced by resolutions with at most two levels,
and the conclusion of  theorem 1.9 follows from theorem 1.5.   

\smallskip
If there is no maximal complexity resolution in the higher rank MR diagram, that has a weak test sequence $\{\varphi_s\}$, such that for every automorphism
$\tau \in Aut(G)$ all the convergent pairs of subsequences, $\{\varphi_s\}$ and $\{\varphi_s \circ \tau\}$, converge into maps:
$\nu_{\tau}:RedComp(CRes_1) \to RedFCl(CRes)$ that satisfy part (5) in lemma 1.13, we gradually replace the higher rank MR diagram with a new diagram in which there
is such a maximal complexity resolution.

We do that by either adding (finitely many) maximal complexity resolutions   that have weak test sequences so that composition with no automorphism
$\tau \in Aut(G)$ can reduce their complexity using lemma 1.13 (i.e., for which part (5) in the lemma holds),  or by
finding a new higher rank MR diagram with resolutions of strictly smaller complexities than the maximal complexity of resolutions in the previous diagram.
Since every descending chain of complexities of resolutions terminates after finitely many steps, this replacement procedure terminates after finitely
many steps. When it terminates, there will exist maximal complexity resolutions that satisfy the condition on the existence of weak test sequences
for which for every automorphism $\tau \in Aut(G)$, and every convergent subsequences, the obtained maps $\nu_{\tau}$ satisfy property (5) in lemma 1.13.
  
We look at all the possible sequences of automorphisms $\{\varphi_s\}$ from $Aut(G)$. From each such sequence we can pass to a convergent subsequence,
that converges into some $m$-collection with a single non-trivial resolution, $Res$. 

Let $Compl_{max}$ be the maximal
complexity of the non-trivial resolutions in the original higher rank MR diagram. If the complexity of $Res$ is strictly bigger than $Compl_{max}$,
then the argument that was used to prove lemma 1.13, implies that conclusions (1)-(4) in lemma 1.13 imply for the resolution $Res$.
Hence, after possibly passing to a subsequence, and iteratively reducing the complexities of the obtained resolutions, 
it is possible to replace the resolution $Res$ with a resolution $Res_1$, such that $Res_1$ has complexity bounded by $Compl_{max}$.

Suppose that the complexity of $Res_1$ is $Compl_{max}$. Once again we apply the argument that was used in the proof of lemma 1.13. If one of the
conclusions (1)-(4) apply to either a subsequence of $\{\varphi_s\}$ or for a subsequence of $\{\varphi_s \circ \tau\}$ for some fixed automorphism
$\tau \in Aut(G)$, then it is possible to replace the resolution $Res_1$ with a resolution of strictly smaller complexity. Otherwise a subsequence of
the automorphisms $\{\varphi_s\}$ is a weak test sequence of $Res_1$, and for every automorphism $\tau \in Aut(G)$, it is possible to obtain
an automorphism $\nu_{\tau}$ that satisfies conclusion (5) in lemma 1.13.

Therefore, if it was not possible to find a resolution $Res$ with a complexity $Compl_{max}$, and with a weak test sequence of automorphisms,
$\{\varphi_s\}$, such that for every automorphism $\tau \in Aut(G)$, and every pair of convergent sequences
it is possible to construct a map $\nu_{\tau}$ that satisfies conclusion (5)
in lemma 1.13, then from a sequence of automorphisms, $\{\varphi_s\}$, it is possible to pass to a subsequence which is a weak test sequence of some resolution
that has complexity that is strictly smaller than $Compl_{max}$.

At this point we can reconstruct a new higher rank diagram, repeating the construction in theorem 2.8 in [Se8]. We go over all the $m$-collections
of cover resolutions that have a weak test sequence of automorphisms, and for which their only non-trivial resolution has complexity that is strictly
bounded by $Compl_{max}$. Clearly, there are only countably many such $m$-collections. 
By the compactness argument that was used in the proof of theorem 2.8 in [Se8], finitely many such $m$-collections of cover resolutions suffice. Hence, we 
constructed a new higher rank MR diagram that contain finitely many $m$-collections, and the only non-trivial cover resolution in each $m$-collection has
complexity that is strictly smaller than $Compl_{max}$.

Starting with the new diagram we repeat the whole construction. Either we find an $m$-collection of covers in which the only non-trivial resolution has
complexity that is equal to the maximal complexity of the resolutions in the new diagram, and this resolution has a weak
test sequence, such that for every automorphism $\tau \in Aut(G)$, the corresponding map $\nu_{\tau}$ satisfies conclusion (5) in lemma 1.13, or we
can construct a new higher rank MR diagram with strictly smaller maximal complexity resolutions.

The descending chain condition for complexities of resolutions imply that this procedure terminates after finitely many steps, and our previous argument imply that
when it terminates the maximal complexity resolutions have only two levels (in their reduced modeled resolutions)
In this case the conclusions of theorem 1.9 follow from theorem 1.5.

\line{\hss$\qed$}

So far we have assumed that in each $m$-collection in the higher rank MR diagram there is at most one non-trivial resolution, 
and that the virtually abelian decompositions along non-trivial resolutions contain no edges with finite (nor trivial) edges groups.
We further assumed that the quotient maps
along the non-trivial resolutions are proper quotients, which guarantee that each resolution (or $m$-collection) in the MR diagram has a weak test
sequence of automorphisms that factor through it.

The next theorem drops the assumption for proper quotient maps. It still requires that the virtually abelian decompositions along the non-trivial resolutions 
contain no edges with finite (nor trivial) edge groups. This assumption is kept in this paper, and will be dropped only in the next paper in this sequence. 

\proclaim{Theorem 1.15 (cf. theorem 1.9)} With the assumptions of lemma 1.1, suppose that   
in each of the  (finitely many) $m$-collections of  cover resolutions in the higher rank Makanin-Razborov diagram, there is at most one resolution
that has more than a single level, and
all the virtually abelian decompositions that are associated with the various levels 
of the resolutions that have more than a single level do not contain edges with finite (nor trivial) edge groups.

Then conclusions (1) and (2) in theorem 1.8 hold. In particular, $H$ admits a higher rank JSJ decomposition.
\endproclaim

\nfp In theorem 1.9 we assumed that in all the resolutions in the higher rank MR diagram, all the quotient maps along the non-trivial resolutions are
proper quotients. By the construction of the higher rank MR diagram [Se8], this assumption guarantees that each resolution (or rather an $m$-collection) in the 
higher rank MR diagram has a weak test sequence of automorphisms that factor through it. 

Since weak test sequences are used to prove theorem 1.9, to apply the argument that is used in the proof of theorem 1.9 to prove theorem 1.15, we need 
to modify the construction of the higher rank MR diagram that appears in section 3 of [Se8], to guarantee that resolutions, or $m$-collections, that appear in the higher
rank diagram, have weak test sequences, even if  some of the quotient maps along the resolutions in an $m$-collection are not proper quotients.

The construction of resolutions in the higher rank MR diagram in [Se8] guarantees that resolutions, and hence $m$-collections, have weak test sequences
that factor through them, only if all the quotient maps along the resolutions are proper quotients.
Therefore, to guarantee that all the $m$-collections in the higher rank MR diagram have weak test sequences that factor through them, we need to modify the
construction of the resolutions in the $m$-collections in the higher rank MR diagram. Indeed, a construction of resolutions that admit weak test sequences even when
some of the quotient maps along the resolutions are isomorphisms (and not proper quotients), appears in our work on Makanin-Razborov diagrams of pairs,
that analyze and encode varieties over a free semigroup [Se9].

The construction of resolutions that is used in [Se9], that involves infinite chains of decompositions that are eventually replaced by a finite resolution,
construct finite resolutions with weak test sequences, where some of the quotient maps between consecutive levels are isomorphisms and some are
proper quotients. 

When we apply the construction of resolutions that is used in [Se9] (for pairs and semigroups), we start with a sequence of automorphisms, $\{\varphi_s\}$
in $Aut(G)$, and pass to a subsequence which is a weak test sequence w.r.t.\ some finite resolution of a limit quotient of $H$ (the characteristic
finite index subgroup in $G$). 

The same compactness argument that is used in the construction of the higher rank MR diagram in section 2 in [Se8], implies that using the construction of resolutions
in [Se9], there exists a higher rank MR diagram, that contains finitely many $m$-collections of  resolutions, and for each $m$-collection there is
a weak test sequence of automorphisms that factor through it.

Once we replaced the higher rank MR diagram with a diagram in which every $m$-collection has a weak test sequence, the proof of theorem 1.15 follows by
the same argument that was used in the proof of theorem 1.9.

\line{\hss$\qed$}

\vglue 1.5pc
\centerline{\bf{\S2. A higher  rank JSJ decomposition of a product II:}} 
\centerline{\bf{general diagrams}}
\medskip
\medskip

In the previous section we analyzed higher rank MR diagrams of products of hyperbolic spaces, 
in which every $m$-collection contains at most a single resolution with more than 2 
levels. In this case we were able to associate with the product a higher rank JSJ decomposition,
from which it is possible to extract information on the dynamics of individual automorphisms, and on the algebraic
structure of $Out(G)$.

In this section we continue with products of hyperbolic spaces, but we analyze general higher rank MR diagrams. 
In this case we do not manage to get  graphs of groups decompositions (of limit quotients) that are associated with the various factors,
and are invariant under a finite index subgroup of $Out(G)$. The techniques that we use for analyzing general higher rank MR diagrams of products,
does not give us an invariant graph of groups, but we do construct (non-invariant) graphs of groups that enable us to analyze the dynamics of individual
automorphisms and the algebraic structure of $Out(G)$ as in the previous section.


We continue with the notation of the previous section. Hence, $X$ is a product of $m$ hyperbolic spaces $V_1,\ldots,V_m$, 
and $G$ is an HHG that acts on $X$ properly 
and cocompactly. $H$ is a characteristic finite index subgroup of $G$ that preserves the factors $V_1,\ldots,V_m$ and acts on them weakly acylindrically and
isometrically. 

Following section 2 in the previous paper in this sequence [Se8], with the action of $H$ on $X$ we associate a higher rank MR diagram (theorem 2.8 in [Se8]). 
The higher rank diagram 
contains finitely many $m$-collections of cover resolutions, and each automorphism in $Aut(G)$ factors through at least one of the $m$-collections in
the higher rank MR diagram. As in the previous section we analyze general higher rank MR diagrams of products 
by gradually relaxing the assumptions on their structure.

\vglue 1.5pc
\proclaim{Theorem 2.1} Suppose that   
in each of the  the (finitely many) collections of $m$ cover resolutions in the higher rank Makanin-Razborov diagram, each resolution has at most two levels.

Definition 1.6 associates a complexity with a 2-level resolution. Given an $m$ collection in which all the resolutions have at most two levels,
we set its complexity to be the sequence of complexities of its $m$ resolutions, ordered in a non-increasing order. The set of complexities of  $m$-collections
can be ordered by ordering the complexities in a lexicographical order. 

Since we assumed that all the resolutions in the $m$-collections in the higher rank MR diagram have at most two levels, there are maximal complexity 
$m$-collections in the higher rank MR diagram. 

With $H$ and the QH and the virtually abelian groups in the virtually cyclic decompositions that are associated with
the $m$ resolutions in the maximal complexity $m$-collections we associate a groupoid $GD1$.
The objects in $GD1$ are all the possible ordered collections of QH and virtually abelian vertex groups in the maximal complexity $m$-collections, where QH vertex 
groups are equipped with markings for their boundaries. Hence, the groupoid $GD1$ has finitely many objects.

With each outer automorphism $\tau \in Out(G)$ we associate finitely many maps (arrows)  $\mu^j_{\tau}$ between the QH and the virtually abelian vertex groups 
that are associated with objects (vertices) in $GD1$.   
\roster
\item"{(i)}"  A map $\mu^j_{\tau}$ consists of isomorphisms from  the QH and the virtually abelian vertex groups that appear in the virtually
abelian decompositions that are associated with the collection in the source vertex 
in $GD1$ to the QH and virtually abelian vertex groups that are associated with the target vertex in $GD1$. The isomorphisms are defined up to compositions 
with inner automorphisms.

\item"{(ii)}" With each arrow $\mu^j_{\tau}$ in $GD1$ there is an inverse arrow $\mu^j_{\tau^{-1}}$.

\item"{(iii)}" if $\tau=\tau_1 \circ \tau_2$, then for each arrow $\mu^j_{\tau}$, there exist arrows $\mu{j_1}_{\tau_1}$ and $\mu^{j_2}_{\tau_2}$ such
that: $\mu^j_{\tau}=\mu^{j_1}_{\tau_1} \circ \mu^{j_2}_{\tau_2}$.
\endroster
\endproclaim

\nfp 
Let $Coll$ be one of the  maximal complexity $m$-collections in the higher rank MR diagram.
Since the resolutions in the $m$-collection have at most two levels, they have a weak test sequence. 
Let $\{ \varphi_s \}$ be a weak test sequence w.r.t.\ the $m$-collection $Coll$. Note that this implies that the sequence $\{ \varphi_s \}$ restricts to
weak test sequences of all the $m$ resolutions in the collection $Coll$. 

Let $\tau \in Aut(G)$.
Since $\{ \varphi_s \}$ is a weak test sequence w.r.t.\ a maximal complexity collection, the sequence $\{ \varphi_s \circ {\tau}^{-1} \}$ must contain a subsequence that 
is a weak test sequence of at least one of the maximal complexity collections in the higher rank MR diagram.

\vglue 1.5pc
\proclaim{Proposition 2.2} Let $\{\varphi_s\}$ be a weak test sequence of a maximal complexity 
$m$-collection $Coll$, let $\tau \in Aut(G)$, and suppose that $\{ \varphi_s \circ {\tau}^{-1}\}$
factors through a maximal complexity $m$-collection $Coll_1$ from the higher rank MR diagram. Then:

\roster
\item"{(1)}" $Coll_1$ is a maximal complexity $m$-collection.

\item"{(2)}" the sequence $\{ \varphi_s \circ {\tau}^{-1}\}$ 
contains a subsequence that 
is a weak test sequence w.r.t.\  $Coll_1$.

\item"{(3)}" the convergent subsequence of the sequence $\{ \varphi_s \circ {\tau}^{-1}\}$, which is a weak test sequence of the collection $Coll_1$, defines  
maps that we denote  $\mu_{\tau}$ from the QH vertex groups and the virtually abelian vertex groups in $Coll$  onto the QH vertex groups and virtually abelian vertex
groups in $Coll_1$. These maps are defined up to a composition with inner automorphisms.
\endroster
\endproclaim

\nfp The weak test sequence, $\{\varphi_s\}$, converges into limit groups with virtual cyclic decompositions that are associated with
the maximal complexity $m$-collection $Coll$. 

Let $\tau \in Aut(G)$.
A subsequence of $\{\varphi_s \circ {\tau}^{-1}\}$ factors through some $m$-collection $Coll_1$. We still denote the subsequence,
$\{ \varphi_s \circ {\tau}^{-1}\}$. Since $\tau \in Aut(G)$, $\{ \varphi_s \circ {\tau}^{-1}\}$ converges
into limit groups with virtually cyclic decompositions that have the same structure as the virtually cyclic decompositions that are associated with the
$m$-collection $Coll$. Since the complexities of of these  virtually cyclic decompositions are bounded by the complexities of the virtually cyclic
decompositions that are associated with $Coll_1$, and since they have the complexities of the maximal complexity collection $Coll$, $Coll_1$ must be
a maximal complexity collection as well. Furthermore, the virtually cyclic decompositions that are associated with the limit of the
sequence, $\{ \varphi_s \circ {\tau}^{-1}\}$, must be the virtually cyclic decompositions that are associated with the $m$-collection $Coll_1$. Hence,
$\{\varphi_s \circ {\tau}^{-1}\}$ is a weak test sequence for the $m$-collection $Coll_1$.

Since $\{\varphi_s\}$ and $\{ \varphi_s \circ {\tau}^{-1}\}$ converge into limit groups with virtually abelian groups that have the same structure, the automorphism
$\tau$ naturally maps the QH and virtually abelian vertex groups in the virtually cyclic decompositions that are associated with
$\{\varphi_s\}$ onto the QH and virtually abelian vertex groups in the virtually cyclic decompositions that are associated with 
$\{\varphi_s \circ {\tau}^{-1}\}$.
Since QH and virtually abelian groups have the Hopf property, the homomorphisms that are associated with $\tau$, that are necessarily 
epimorphisms, must be isomorphisms.

\line{\hss$\qed$}

Proposition 2.2 proves that from any weak test sequence $\{ \varphi_s \}$ w.r.t.\ a maximal complexity collection $Coll$, and an automorphism $\tau \in Aut(G)$, it is possible
to pass to a subsequence for which the precomposed sequence $\{ \varphi_s \circ {\tau}^{-1} \}$ is a weak test sequence w.r.t.\ some maximal complexity collection
$Coll_1$. the two convergent sequences $\{ \varphi_s \}$ and $\{\varphi_s \circ {\tau}^{-1} \}$, enable us to define a map $\mu_{\tau}$ between the QH and 
virtually abelian 
vertex groups in the virtually cyclic decompositions that are associated with the $m$ resolution in $Coll$, to those in $Coll_1$, where the map is well-defined up
to composition with an inner automorphism. 

The next proposition proves that for each $\tau \in Aut(G)$ there can only be finitely many distinct maps $\mu_{\tau}$ that are constructed from  all the weak test sequences
w.r.t.\ all the maximal complexity collections in the higher rank MR diagram. 

\vglue 1.5pc
\proclaim{Proposition 2.3} Let $\tau \in Aut(G)$ and let $Coll_1$ and $Coll_2$ be two maximal complexity $m$-collections in the higher rank
MR diagram.  Up to composition with inner automorphisms, there  exist only finitely many maps $\mu_{\tau}$, 
that are constructed according to part (iii) of proposition 2.2
from  weak test sequences, $\{ \varphi_s \}$ w.r.t.\ $Coll_1$ and  $\{ \varphi_s \circ {\tau}^{-1} \}$ w.r.t.\ $Coll_2$.
\endproclaim

\nfp To prove the proposition we use a compactness argument that is similar to the one that was used in the construction of the higher rank
MR diagram.

Let $\tau \in Aut(G)$.
We look at the collection of all the sequences, $\{ \varphi_s \}$, that factor through $Coll_1$ and converge into  limit groups with  virtually cyclic
decompositions that have the same structure as the virtually cyclic decompositions that are associated with $Coll_1$.
Furthermore, we require that  $\{ \varphi_s \circ {\tau}^{-1} \}$ factor through $Coll_2$ and converge to limit groups with virtually cyclic decompositions that 
have the same structure as the virtually cyclic decompositions that are associated with $Coll_2$. Note that we don't require these sequences to be weak
test sequences w.r.t.\ $Coll_1$ and $Coll_2$. 

By the same argument that was used to deduce part (iii) of proposition 2.2 (for weak test sequences), each pair of such sequences gives  
maps from the collection of QH and virtually abelian vertex groups
in the virtually cyclic decompositions that are associated with $Coll_1$ isomorphically onto the collection of QH and virtually abelian vertex groups in the 
virtually cyclic decompositions that are associated with $Coll_2$. 

If there are only finitely many such maps, the conclusion of proposition 2.3 follows. Hence, we assume that there are infinitely many such maps (that are defined
up to compositions with inner automorphisms). Clearly there can  only be countably many such maps, and we order them.

We start with the given sequence of sequences, $\{ \varphi^1_s \}, \, \{varphi^2_s\}, \ldots$, and the maps $\{ \mu^j_{\tau} \}$ between the QH and virtually abelian
vertex groups in $Coll_1$ and $Coll_2$ that are associated with them. By our assumptions, the maps $\mu^j_{\tau}$ are distinct even after 
compositions with inner automorphisms.

We construct a new sequence. For each index $j$, we choose an automorphism $\psi_j \in Aut(G)$, such that:

\roster
\item"{(1)}"  $\psi_j$ factors through $Coll_1$ and $\psi_j \circ {\tau}^{-1}$ factors through $Coll_2$.

\item"{(2)}" the pair $\psi_j$ and $\psi_j \circ {\tau}^{-1}$ do not factor through the maps $\mu^i_{\tau}$, for $i=1,\ldots,j$.   

\item"{(3)}" the projections to the various factors of the images of $H$ under $\psi_j$ and $\psi_j \circ {\tau}^{-1}$ 
satisfy the relations and the inequalities in the (cover) limit groups that
are associated with $Coll_1$ and $Coll_2$ for all the elements in a ball of radius
$j$ in $H$ (w.r.t.\ a fixed set of generators of $H$). Furthermore, elements in the ball of radius $j$ in $H$ that are not mapped 
to the bounded vertex groups in the projections to the various factors in $Coll_1$ and $Coll_2$, are mapped by $\psi_j$ and $\psi_j \circ {\tau}^{-1}$ 
to elements of length at least $j$ 
times the bound on the lengths of the images of the (fixed) generators of the bounded vertex groups in $Coll_1$ and $Coll_2$ for each of the factor spaces  of $X$.   
\endroster

Such automorphisms $\psi_j$ exist for every index $j$, since we assumed that there are infinitely many distinct maps $\mu^j_{\tau}$, so we can 
pick $\psi_j$ to be an automorphism from the sequence: $\{\varphi^{j+1}_s\}$, for large enough $s$.

From the automorphisms, $\{ \psi_j \}$, we can pass to a convergent subsequence. This subsequence converges into limit groups with virtually cyclic 
decompositions that have
the same structure as the virtually cyclic decompositions that are associated with the maximal complexity $m$-collection $Coll_1$. Hence, with the 
convergent sequences (still denoted), $\{\psi_j\}$ and $\{\psi_j \circ {\tau}^{-1}\}$, it is possible to associate  maps $\nu_{\tau}$ between the QH and
virtually abelian vertex groups in $Coll_1$ to these in $Coll_2$. But this map must be one of the maps $\mu^{j_0}_{\tau}$ that were associated with $\tau$
in the ordered sequence. Therefore, for large $j>j_0$, the automorphisms $\psi_j$ and $\psi_j \circ \tau$, factor through $\mu^{j_0}_{\tau}$, a contradiction to
the way the automorphisms $\psi_j$ were chosen. Hence, there can only be finitely many maps $\mu^j_{\tau}$ and the proposition follows.

\line{\hss$\qed$}

Propositions 2.2 and 2.3 enable us to complete the proof of theorem 2.1. The finitely many objects in the groupoid $GD1$, 
are associated with all the possibilities for 
collections of QH vertex groups (with marked boundary components) and virtually abelian vertex groups in all the virtually cyclic decompositions that
are associated with resolutions in the maximal complexity $m$-collections in the higher rank MR diagram.

With each automorphism $\tau \in Aut(G)$ we associate the finitely many maps $\mu^j_{\tau}$, that are constructed from weak test sequences
$\{ \varphi_s \}$ and $\{ \varphi_s \circ {\tau}^{-1} \}$  w.r.t.\ maximal complexity collections $Coll_1$ and $Coll_2$. Proposition 2.2 explains how to construct
these maps, and proposition 2.3 proves that there are only finitely many of them for each $\tau \in Aut(G)$. By construction, the maps 
$\mu^j_{\tau}$ depend only on the class of $\tau$ in $Out(G)$.

By construction, for every map $\mu^j_{\tau}$ there is a map $\mu^j_{\tau^{-1}}$, so that their compositions in both orders are the corresponding identity maps.
Furthermore, if $\tau = \tau_1 \circ \tau _2$, then for each of the maps $\mu^j_{\tau}$, there exist maps, $\mu^{j_1}_{\tau_1}$ and $\mu^{j_2}_{\tau_2}$, such
that:   
$\mu^j_{\tau}=\mu^{j_1}_{\tau_1} \circ \mu^{j_2}_{\tau_2}$. This concludes the proof of theorem 2.1

\line{\hss$\qed$}

In case all the resolutions in the higher rank MR diagram have at most two levels, theorem 2.1 associates with the maximal complexity collections in the
diagram a groupoid. The groupoid is non-trivial if there are QH or virtually abelian vertex groups in the virtually cyclic decompositions that are associated
with maximal complexity collections in the diagram. The next theorem associates another groupoid with higher rank MR diagrams, for which the virtually
cyclic decompositions that are associated with their  maximal complexity
collections contain no QH nor virtually abelian vertex groups.  

\vglue 1.5pc
\proclaim{Theorem 2.4} With the assumptions of theorem 2.1, suppose that the virtually cyclic decompositions that are associated with maximal complexity
collections in the higher rank MR diagram contain no QH nor virtually abelian vertex groups.

With $H$ and the infinite virtually cyclic edge groups in the virtually cyclic decompositions that are associated with
the $m$ resolutions in the maximal complexity $m$-collections we associate a groupoid $GD2$.
The objects in $GD2$ are  the ordered collections of virtually infinite cyclic edge groups in the maximal complexity $m$-collections in the higher rank MR diagram.


Suppose that there are $n$ edges with virtually cyclic edge groups in the maximal complexity $m$-collections in the higher rank MR diagram.
With each outer automorphism $\tau \in Out(G)$ we associate finitely many arrows, that we denote:  $\nu^j_{\tau}$, between pairs of vertices in $GD2$. 
With each arrow, $\nu^j_{\tau}$, we  associate an element in $Z^n$ 
that we denote
$\ell^j_{\tau}$.
    
\roster
\item"{(i)}" With each arrow $\nu^j_{\tau}$ in $GD2$ there is an inverse arrow $\nu^j_{\tau^{-1}}$. $\ell^j_{\tau}+\ell^j_{\tau^{-1}}$ is a uniformly bounded
element in $Z^n$. 

\item"{(ii)}" if $\tau=\tau_1 \circ \tau_2$, then for each arrow $\nu^j_{\tau}$, there exist arrows $\nu{j_1}_{\tau_1}$ and $\nu^{j_2}_{\tau_2}$ such
that the composition of the arrows that are associated with $\nu^{j_1}_{\tau_1}$ and $\nu^{j_2}_{\tau_2}$ is the arrow $\nu^j_{\tau}$,
and $\ell^j_{\tau}-\ell^{j_1}_{\tau_1}-\ell^{j_2}_{\tau_2}$ is a uniformly bounded element in $Z^n$.
\endroster
\endproclaim

\nfp The argument is similar to the proof of theorem 2.1. 
We define the complexity of a virtually abelian decomposition with no QH nor virtually abelian vertex groups, and only two levels, to be the following tuples
of numbers:

\roster
\item"{(1)}" the number of edges with trivial or finite edge groups.

\item"{(2)}" the orders of the finite edge groups, from the smallest order to the biggest one (counted with multiplicities).

\item"{(3)}" the number of edges with infinite virtually cyclic edge groups.
\endroster 

We order the complexities lexicographically, a virtually abelian decomposition $\Lambda_1$ has bigger complexity than $\Lambda_2$, if the number of edges
in (1) in $\Lambda_1$ is bigger, or if there is equality in (1)  and the tuples of orders in (2) of $\Lambda_1$ is smaller in lexicographical order than
that of $\Lambda_2$, or if the tuples in (1) and (2) for $\Lambda_1$ and $\Lambda_2$ are equal, and the number of edges in (3) in $\Lambda_1$ is bigger
than that of $\Lambda_2$.

We set the complexity of an $m$-collection of resolutions with at most two levels, with no QH nor virtually abelian vertex groups, to be the $m$-tuple
of complexities of the virtually abelian decompositions that are associated with the $m$ (cover) resolutions in the $m$-collection ordered in a non-increasing order.
The complexities of such $m$-collections are naturally ordered lexicographically. 

Given a higher rank MR diagram where all the resolutions in the $m$-collections have at most two levels, and all the virtually abelian decompositions that are
associated with the resolutions have no QH nor virtually abelian vertex groups, we look at the maximal complexity $m$-collections in the higher rank diagram.
If $\{\varphi_s\}$ is a weak test sequence of such a maximal complexity $m$-collection, and $\tau \in Aut(G)$, then $\{\varphi_s \circ {\tau}^{-1}\}$ must have a subsequence
that is a weak test sequence of some maximal complexity $m$-collection in the higher rank diagram as well.

Given two maximal $m$-collections $Coll_1$ and $Coll_2$ from the higher rank MR diagram, and $\tau \in Aut(G)$, 
We look at all the weak test sequences, $\{\varphi_s\}$ of $Coll_1$, for which $\{\varphi_s \circ {\tau}^{-1}\}$ is a weak test sequence of $Coll_2$. 

A weak test sequence
$\{\varphi_s\}$ of $Coll_1$ converges into limit groups with the virtually cyclic decompositions that are associated with $Coll_1$. Similarly 
$\{\varphi_s \circ {\tau}^{-1}\}$
converges into limit groups with the virtually cyclic decompositions that are associated with $Coll_2$. Hence, $\tau$ maps the vertex and edge groups in the virtually cyclic
decompositions 
that are associated with $Coll_1$ onto the vertex and edge groups in the virtually cyclic decompositions that are associated with $Coll_2$.
We denote this map $\nu_{\tau}$.

With each $m$-collection, $Coll$, we fix (finite) generating sets of all the vertex groups, and a generator of the maximal cyclic subgroup in all the
 edge groups in the $m$ virtually cyclic decompositions that are associated with the $m$-collection.
By construction, each automorphism $\varphi \in Aut(G)$ that factors through an $m$-collection $Coll$, maps the fixed set of generators of a vertex group 
into elements, that map some point in the corresponding projection space, a uniformly bounded distance. i.e., a distance that is bounded by a constant
that depends only on the fixed
generating set and not on the specific automorphism $\varphi$.

Let $\{\varphi_s\}$ and $\{\varphi_s \circ {\tau}^{-1}\}$, be two convergent weak test sequences that factor through the maximal complexity
$m$-collections $Coll_1$ and $Coll_2$. Then the vertex and edge groups in the virtually cyclic decompositions that are associated with the limits of 
$\{\varphi_s \}$ are mapped by a map that we denote $\nu_{\tau}$ to the vertex and edge groups in the virtually cyclic decompositions that
are associated with the limits of $\{\varphi_s \circ {\tau}^{-1}\}$.

Given an edge $E$ in one of the virtually abelian decompositions in the $m$-collection $Coll_1$, where $E$ has an infinite virtually cyclic edge group,
we look at the points $p^1_s,p^2_s$ 
that move minimally by the images 
under $\{\varphi_s\}$ of the
fixed generating sets of the vertex groups that are adjacent to the edge, and at the points $q^1_s,q^2_s$ 
that move minimally by the images under $\{\varphi_s \circ {\tau}^{-1}\}$
of the vertex groups in the virtually cyclic decompositions in the $m$-collection $Coll_2$ that are mapped by $\nu_{\tau}$ to these two vertex groups.  

By possibly composing the automorphism $\tau$ with an inner automorphisms, we can assume that the distance between $p^1_s$ and $q^1_s$ is uniformly bounded.
Furthermore, there is some power $\ell_E$, such that $\varphi_s ({g_E}^{\ell_E})(p^2_s)$ has a uniformly bounded distance from $q^2_s$, where $g_E$ is a fixed generator of
the maximal cyclic group in the virtually cyclic edge group that stabilizes the edge $E$. Note that the power $\ell_E$ is not uniquely defined, but it is 
defined up to some universal constant  (that does not depend on the automorphism $\tau$).

Hence, with the convergent weak test sequences, $\{\varphi_s\}$ and $\{\varphi_s \circ {\tau}^{-1}\}$, it is possible to associate a tuple of integers with the edges
in the virtually cyclic decompositions that are associated with the $m$-collection $Coll_1$. Therefore, with these two weak test sequences it is possible to
associate an element in $Z^n$, where $n$ is the total number of edges with infinite virtually cyclic edge groups in the virtually abelian decompositions
that are associated with a maximal complexity collection in the higher rank MR diagram. We denote this element in $Z^n$, $\ell_{\tau}$. Note that
$\ell_{\tau}$ is not uniquely defined, but it is defined up to a uniformly bounded element in $Z^n$. 
  
Properties (i) and (ii) of the elements $\ell^j_{\tau}$ that are associated with an automorphism $\tau \in Aut(G)$ and all the possible convergent weak
test sequences, $\{\varphi_s\}$ and $\{\varphi_)s \circ {\tau}^{-1}\}$, follow from the construction of the elements $\ell^j_{\tau}$. It is left to prove that with
$\tau$ it is possible to associate only finitely many arrows, i.e., finitely many such elements $\ell^j_{\tau}$ (for all the convergent weak test sequences). 

Suppose that there is an automorphism $\tau \in Aut(G)$, with an infinite sequence of elements $\ell^j_{\tau}$. In this case we can pass to an unbounded
subsequence of elements $\ell^j_{\tau}$ in $Z^n$.

We look at pairs of convergent weak test sequences of maximal complexity $m$-collections $Coll_1$ and $Coll_2$:  $\{\varphi^j_s\}$ and 
$\{\varphi^j_s \circ {\tau}^{-1}\}$, with the unbounded sequence of elements
$\ell^j_{\tau} \in Z^n$.  
From each sequence $\{\varphi^j_s\}$ we pick an automorphism $\psi_j$, such that the sequences, $\{\psi_j\}$ and $\{\psi_j \circ {\tau}^{-1}\}$, converge.
Since the sequence of elements $\{\ell^j_{\tau}\}$ is unbounded, we can pass to a further subsequence, for which at least one of 
the virtually cyclic decompositions 
that is associated
with the pair of sequences is a proper refinement of the corresponding virtually cyclic decomposition that is associated with the 
$m$-collections $Coll_1$ and $Coll_2$. Since these
$m$-collections are of maximal complexity, we got a contradiction, since the $m$-collection of virtually cyclic decompositions that are associated with
the convergent subsequences, $\{\psi_j\}$ and $\{\psi_j \circ {\tau}^{-1}\}$, have strictly bigger complexity. Therefore, with each automorphism $\tau \in Aut(G)$
there are at most finitely many associated arrows $\nu^j_{\tau}$ in the groupoid $GD2$.

\line{\hss$\qed$}

Theorems 2.1 and 2.4 analyze higher rank MR diagrams in which all the resolutions in the $m$-collections in the diagram have no more than two levels. In both
theorems a higher rank JSJ is constructed. However,  we only managed to associate with the outer automorphism group, $Out(G)$, and with the higher
rank JSJ, two groupoids, $GD1$ and $GD2$. This is mainly since we were not able to associate with elements in $Out(G)$ unique maps or (outer) automorphisms
of the QH and virtually abelian vertex groups that appear in the higher rank JSJ, but only finitely many outer automorphisms, and we were not able to associate 
unique Dehn twists with elements in $Out(G)$ and the various edges with virtually cyclic edge groups in the higher rank JSJ, but only finitely many ones.

In the sequel we do manage to get uniqueness of the maps, not into the mapping class groups of the punctured orbifold, but rather to orbifolds that are obtained by attaching disks 
to the boundaries (or some s.c.c.) of the 2-orbifolds that are associated with the QH vertex groups in the higher rank JSJ decomposition. Combining this with theorem
2.1, this means that 
for each automorphism $\tau \in Aut(G)$, these unique maps to the mapping class groups of the closed 2-orbifolds lift to at most finitely many maps of the mapping class
groups of the punctured orbifolds.

\medskip
We continue by further applying
 the arguments that were used in the previous section to generalize the results to $m$-collections with resolutions with arbitrary many levels.
In this paper we still assume that all the resolutions that appear in $m$-collections in the higher rank MR diagram
(that was constructed in theorem 2.8 in [Se8]) are of minimal rank.
This means that all the virtually abelian decompositions along these 
resolutions have no edges with finite (nor trivial) edge groups.
General $m$-collections and their higher rank JSJ decomposition  will be analyzed and constructed in the next paper in this sequence.

\proclaim{Theorem 2.5 (cf. theorems 1.9 and 1.15)} Suppose that  in every $m$-collection in the higher rank MR diagram of a group $G$ that acts on a product
of hyperbolic spaces $X$ (where the action satisfies the properties that appear in lemma 1.1),   
all the virtually abelian decompositions that appear along the $m$ cover resolutions do not have edges with finite 
nor trivial edge groups
(i.e., the cover resolutions, and hence the $m$-collection, are of minimal rank).
 

Then it is possible to associate with the higher rank diagram a higher rank JSJ decomposition, i.e., a virtually abelian (JSJ like) decomposition of a quotient of $H$ that is
associated with each factor, and  groupoid, similar to either the groupoid $GD1$ that was constructed in theorem 2.1,
or $GD2$ that was constructed in theorem 2.4.
\endproclaim

\nfp We start by assuming that all the quotient maps along the resolutions in the higher rank MR diagrams are proper quotients. In this case all the
$m$-collections in the diagram have weak test sequences. 
We follow the analysis of resolutions that we used in the previous section. i.e., in case each $m$-collection has at most a single
non-trivial resolution (theorem 1.9), and apply it to the $m$ resolutions in the $m$-collections simultaneously.



We repeat what we did in the proof of theorem 1.9.
We modify the resolutions in each $m$-collection from the higher rank diagram using the procedure that is applied in the proof of
theorem 1.9, and push QH and virtually abelian vertex groups to lower levels, if the edge groups that are connected to these vertex groups
are all elliptic in the level below them. We further replace the completion of each resolution with its modeled completion (definition 1.10). 

With the  $m$-collections of modeled completions in the higher rank diagram we associated a complexity, which is the $m$-tuple of 
complexities of the modeled completions in the
$m$-tuple in a non-increasing order (the complexity of a modeled completion is the complexity that is used in the proof of theorem 1.9).
The complexities of $m$-collections of modeled completions are well-ordered, and as in theorem 1.9 we continue by analyzing maximal complexity $m$-collections.
  
Let $Coll$ be an  $m$-collection of modeled completions in the (modified) higher rank MR diagram. 
Let $\{\varphi_s\}$ be a weak test sequence w.r.t.\  $Coll$, and let $\tau \in Aut(G)$. By the construction of formal solutions that appears in [Se2],
there exists a  subsequence of pairs:
$\{\varphi_s\}$ and $\{\varphi_s \circ {\tau}\}$, that converge into $m$ homomorphisms: $\nu^i_{\tau}:RedComp(CRes^1_i) \to FCl(CRes_i)$, where
$FCl(CRes_i)$, $i=1,\ldots,m$, are  generalized framed closures of the modeled completions in the $m$-collection $Coll$,
and $RedComp(CRes^1_i)$, $1,\ldots,m$, are the $m$ reduced modeled completions  in some $m$-collection
$Coll_1$ from the modified higher rank MR diagram.

Suppose that for every maximal complexity $m$-collection in the higher rank MR diagram, and every weak test sequence that factors through it,
there exists some automorphism $\tau \in Aut(G)$ (that depends on the sequence), such that the pair of sequences, the weak test sequence and the weak test sequence
twisted by $\tau$, have a convergent subsequence that satisfies one of the conclusions  (1)- (4) in lemma 1.13, for at least one of the maps $\nu^i_{\tau}$, $i=1,\ldots,m$.

In that case we apply the argument that was used to prove theorem 1.9 to the given higher rank MR diagram, 
and replace it by another higher rank diagram with strictly smaller maximal 
complexity $m$-collections. Repeating this procedure iteratively, and using the descending chain condition for complexities of $m$-collections,
we get a new higher rank MR diagram, with smaller maximal complexity $m$-collections,  in which at least one of the maximal complexity $m$-collections 
admits weak test sequences,
such that for every automorphism $\tau \in Aut(G)$, and every convergent subsequence of pairs, part (5) of lemma 1.13 holds. 

\smallskip
We continue with a fixed weak test sequence, $\{\varphi_s\}$, that factor through a maximal complexity $m$-collection $Coll$ in the higher rank MR diagram,
for which for any automorphism $\tau \in Aut(G)$, any convergent subsequence of the sequence $\{\varphi_s \circ {\tau}\}$ factor only through a 
maximal complexity 
$m$-collection (of reduced modeled completions) in the higher rank diagram diagram and satisfies part (5) in lemma 1.13.

Following the proofs of theorem 1.9 and lemma 1.13, given a sequence $\{\varphi_s \circ {\tau}\}$ it is possible to pass to a 
convergent subsequence, $\{\varphi_{s_n} \circ {\tau}\}$, that factors through some maximal
complexity $m$-collection of reduced modeled completions,  $Coll_1$, and there are $m$ associated maps:
$\nu^i_{\tau}:RedComp(Res^i_1) \to FCl(Res^i)$, $i=1,\ldots,m$,  where $Res^i_1$ is the $i$-th resolution in the $m$-collection $Coll_1$,  $Res^i$ is the 
$i$-th resolution
in the $m$-collection $Coll$, and $FCl(Res^i)$ is some (trivially framed) closure of (the modeled completion of) the resolution $Res^i$. 
Furthermore, for large enough index $s$, 
the automorphisms from the convergent subsequence, $\varphi_{s} \circ {\tau}$, extend to values of the $m$-collection (of reduced modeled
completions) $Coll_1$ that factor through the maps 
$\nu^i$, $1 \leq i \leq m$.

We look at a sequences: $\{\varphi_{s_n} \circ {\tau_n}\}$, such that:
\roster
\item"{(1)}"  the indices $\{s_n\}$ are a strictly  increasing sequence.

\item"{(2)}" the automorphisms in the sequence factor
through the reduced modeled completions in a maximal complexity collection $Coll_1$ (that depends on the sequence but not on the index $n$), and do not factor
through $m$-collections that are not of maximal complexity. 

\item"{(3)}" each of the automorphisms in the sequence extends to a specialization that factors through an $n$-th $m$-collection of maps:
$\nu^i_{\tau_n}:RedComp(Res^i_1) \to FCl_n(Res^i)$, $i=1,\ldots,m$,  where $Res^i_1$ is the $i$-th resolution in the $m$-collection $Coll_1$, and $Res^i$ is the 
$i$-th resolution in the $m$-collection $Coll$.
\endroster

We look at the collection of all these sequences. The automorphisms $\varphi_{s_n}$  extend to values of $Comp(Res^i)$, and the automorphisms
$\varphi_{s_n} \circ {\tau_n}\}$ extend to values of $RedComp(Res^i_1)$ $i=1,\ldots,m$. By the argument that was used
to prove theorem 1.9, it follows that up to a permutation of the conjugacy classes of the edge groups in the last completions,  
and elements in the quasi-kernels, 
the virtually cyclic edge groups in the completions, 
$RedComp(Res^i)$ and $RedComp(Res^i_1)$, are mapped to 
conjugate subgroups (for each index $n$).

Hence, we may pass to subsequences for which the permutations  are fixed. We look at subsequences for which the corresponding values
(images) of $Comp(Res^i)$ and $RedComp(Res^i_1)$ converge. 
Each such convergent subsequence converges into a limit group, in which the limit of the values of the completion, $Comp(Res^i)$, which are a weak test sequence,
converge into a closure of that completion. Hence, we analyze such a convergent subsequence using a variation of the analysis of formal limit groups as it
appears in section 2 of [Se2], together with the analysis of closures as it appears in the proof of lemma 1.13.  

Therefore, from each convergent subsequence, $\{\varphi_{s_n} \circ {\tau_n}^{-1}\}$,  we can pass to a further subsequence, that converges into $m$ (formal) resolutions,
where with each factor and each modeled completion, $ModComp(Res^i)$, there is an associated formal resolution. Each formal resolution 
has a similar form
as the formal resolutions that encode formal solutions in section 2 of [Se2], and such a formal resolution terminates in a closure of $ModComp(Res^i)$, that
we denote $FCl(Res^i)$, for  
$i=1,\ldots,m$. 

Furthermore, since the analysis of each formal resolution using the analysis  that appears in the proof of lemma 1.13,
starts with the original reduced modeled completion, $RedComp(Res^i_1)$, i.e., with all the QH and virtually abelian vertex groups in
$RedComp(Res^i_1)$, the complexity of the $i$-th formal
resolution that terminates in $FCl(Res^i)$ is bounded above by the complexity of the resolution $RedComp(Res^i_1)$ 
(that is equal to the complexity of $ModComp(Res^i)$)
 for $i=1,\ldots,m$.
 
If for some index $i$, the complexity of the formal resolution is equal
to the complexity of $RedComp(Res^i_1)$, then the structure of the (modeled completion of the) formal resolution is similar to the structure of the 
reduced modeled completion $RedComp(Res^i_1)$, although the modeled structure may be different (but the different modeled structure does not change
the complexity).
 i.e., the 
formal resolution has QH vertex groups with some modeled subtowers above them, and these QH vertex groups  are similar to those in  $RedComp(Res^i_1)$, 
and their boundaries are conjugate. Furthermore, The virtually abelian vertex groups
in the formal resolution are similar to those in the completion $Comp(Res^i_1)$, and their edge groups are conjugate. 

By the same compactness argument that was used in the construction of the higher rank MR diagram (theorem 2.8 in [Se8]), there exists a finite set of 
$m$-collections of (formal) resolutions, through which all the sequences $\{\varphi_{s_n} \circ {\tau_n}\}$ factor. The complexity of each of the $m$-collections
of (formal) resolutions is bounded above by the maximal complexity of the $m$-collections in the original higher rank MR diagram.

At this point we try to reduce the maximal complexity of the $m$-collections in the new higher rank MR diagram in the same way as what we did in the
proof of theorem 1.9. We look at the maximal complexity $m$-collections and the weak test sequences of the form, $\{\varphi_{s_n} \circ {\tau_n}\}$, that
factor through them. If there exists a maximal complexity $m$-collection with a weak test sequence, such that for every automorphism $\sigma \in Aut(G)$,
every convergent pair of sequences, $\{\varphi_{s_n} \circ {\tau_n}\}$ and $\{\varphi_{s_n} \circ {\tau_n} \circ {\sigma}\}$, satisfies conclusion (5) in lemma
1.13, we continue with this maximal complexity $m$-collection, and this weak test sequence that factors through it. 

If there is no such maximal complexity $m$-collection with such a weak test sequence that factors through it, we use conclusions (1)-(4) in lemma 1.13 to 
strictly reduce 
the complexity of the $m$-collection. We further use the construction that appears in the proof of 1.9 to replace the constructed higher rank MR diagram
with a new higher rank MR diagram, such that the complexities of the  $m$-collections in the new diagram are strictly bounded by the maximal complexity
of the $m$-collections in the previous diagram.

As in the proof of theorem 1.9, a finite process of complexity reductions, enables us to construct a higher rank MR diagram, that has a maximal complexity 
$m$-collection with a weak test sequence,  
$\{\varphi_{s_n} \circ {\tau_n}\}$, that factors through it,
such that for every automorphism $\sigma \in Aut(G)$,
every convergent pair of sequences, $\{\varphi_{s_n} \circ {\tau_n}\}$ and $\{\varphi_{s_n} \circ {\tau_n} \circ {\sigma}\}$, satisfies conclusion (5) in lemma
1.13.

If a maximal complexity $m$-collection in the new diagram has the same complexity as the higher rank diagram that we started with, that has the weak
test sequence $\{\varphi_n\}$, we are done. Otherwise, if a maximal complexity $m$-collection has a strictly smaller complexity than the maximal complexity
of the $m$-collections in the original diagram we continue iteratively. 

We start with a maximal complexity $m$-collection in the new diagram and a fixed weak test sequence of it, $\{\varphi_{s_n} \circ {\tau_n}\}$. We look at all
the sequences, $\{\varphi_{s_{n_r}} \circ {\tau_{n_r}} \circ {\sigma_r}\}$, where $\sigma_r \in Aut(G)$ and the new composed sequences satisfy the same properties
as the ones we imposed on the sequence $\{\varphi_{s_n} \circ {\tau_n}\}$.

In each step of the process we obtain a new higher rank Makanin-Razborov diagram with formal resolutions with an additional layer, and possibly strictly
smaller maximal complexity. Because complexities of $m$-collections satisfy a d.c.c.\ after finitely many steps we get a higher rank Makanin-Razborov
diagram, where each formal resolution in the $m$-collections in the diagram  has finitely many layers, and the diagram has a maximal complexity $m$-collection,
such that all the $m$ modeled (formal) completions in the $m$-collection have two top layers with the same complexities. i.e., the two top layers have 
modeled structures with the same 
QH vertex groups in their bottom levels, and with conjugate boundary components, and similar virtually abelian vertex groups, with conjugate edge groups.
Furthermore the weak test sequence that factor through this maximal complexity collection can not be used to reduce the complexity
by composing it with any fixed automorphism from $Aut(G)$ and apply parts (1)-(4) in lemma 1.13.
  
\smallskip
The $m$ limit groups that are associated with the $m$ factors of the given product space are mapped into the $m$ modeled completions that are
associated with the $m$  formal resolutions in each $m$-collection in the higher rank diagram. Hence, they are also mapped into the reduced modeled completions
that are associated with the $m$ formal resolutions.
 
Given a maximal complexity $m$-collection with the top two layers with the same complexity, we look at the $m$ virtually abelian decompositions 
that the $m$ limit groups inherit from their images into (the top layers of) the $m$ reduced modeled completions in the $m$-collection. There are two possibilities
to define the higher rank JSJ decomposition based on these $m$ decompositions:

\roster
\item "{(1)}" the $m$ virtually abelian decompositions that are inherited by the $m$ limit groups.

\item "{(2)}" the abelian decomposition of subgroups of the $m$ reduced completions, that contain only those QH vertex groups that intersect the images of 
the images of the limit groups that are associated with the factors in conjugates of finite index subgroups. Subgroups of  the virtually
abelian vertex groups in the reduced completions, that intersect contain a subgroup that is generated by conjugates of finite index subgroups that are contained
in the images of the $m$ limit groups. 
\endroster

Note that in option (2) we enlarge the $m$ groups that are associated with the higher rank JSJ, and these groups are not canonical, but the $m$ 
decompositions that we obtain have possibly smaller complexities and they are canonical. Also note that the decompositions in the higher rank JSJ are
virtually f.g.\ abelian, and not virtually cyclic as in the previous section.


With the constructed higher rank JSJ decomposition we can further associate the groupoids $GD1$ or $GD2$ (depends if a virtually abelian decomposition 
contains a QH or a virtually abelian vertex group or not), precisely as we did in case all the resolutions in the $m$-collections have at most two levels,
in theorems 2.1 and 2.4  (note that the objects in the groupoid depend on whether we construct the higher rank JSJ according to option (1) or (2)).   

Along the proof we assumed that the maps along the resolutions in the higher rank MR diagrams that were constructed on the way to construct the
higher rank JSJ decomposition are proper quotients. This guarantees that all the constructed resolutions have weak test sequences
of automorphisms that factor through them.

As in the proof of theorem 1.15, it is possible to drop this assumption by using the construction of resolutions that was used in our
work on solutions to systems of equations over free semigroups [Se9]. The quotient maps along the resolutions that are constructed in this work
need not be proper quotients, still the constructed resolutions have weak test sequences that factor through them. The existence of weak test sequences that 
factor through the constructed resolutions suffice for carrying out all the steps of the proof of theorem 2.5.

\line{\hss$\qed$}

\medskip
So far we have constructed a higher rank JSJ decomposition that is associated  with the automorphism group, $Aut(G)$, of a group $G$ that acts 
on a product of finitely many 
hyperbolic spaces. We further associated two groupoids with $Aut(G)$ and the higher rank JSJ. We associated groupoids and not groups, because we could not
associate unique maps that are associated (canonically) with a given automorphism in $Aut(G)$ and each QH or virtually abelian vertex group
in the higher rank JSJ decomposition, 
and could not associate (canonically)  unique displacements with each automorphism in $Aut(G)$ in case there were no QH nor virtually abelian vertex groups. 
 
Indeed, we did not find a way to guarantee uniqueness in the the groupoids that we constructed, but it is possible to get uniqueness in an associated
groupoid. This will enable us to get a (canonical) homomorphism from a finite index subgroup of $Out(G)$ into the direct sum of mapping class groups
of orbifolds that are associated with the ones that appear in the higher rank JSJ decomposition, and the direct sum of outer automorphism groups
of virtually abelian groups that are associated with the virtually abelian vertex groups in the higher JSJ.

\proclaim{Theorem 2.6 (cf. theorem 2.5)} With the assumptions of theorem 2.5, there exists a finite index subgroup of $Out(G)$, with a canonical homomorphism
into the direct sum of  mapping class groups of 2-orbifolds that are associated with the QH vertex groups in the higher rank JSJ decomposition, and the 
direct sum of general linear groups that are associated with
virtually abelian vertex groups in the higher rank JSJ.
\endproclaim

\nfp To prove theorem 2.6 we construct a groupoid $GD3$, that has finitely many objects, and for each outer automorphism in $Out(G)$ there exists at least
one arrow in the groupoid that is associated with it. Unlike the groupoids $GD1$ and $GD2$, with each outer automorphism  $\tau \in Out(G)$ and any two
objects $A,B$  in the groupoid $GD3$, that are not necessarily distinct,  there is at most a single arrow between the pair of objects $A,B$ that is associated
with $\tau$. 

This uniqueness enables us to look at the fundamental group of the groupoid, which is a finite index subgroup of $Out(G)$ (it is defined up to a composition
with an inner automorphism), and
from the construction of the groupoid it follows that this finite index subgroup of $Out(G)$ has a homomorphism into the direct sum of mapping class groups and
general linear groups as the theorem claims. 

Given a QH vertex group $Q$ in the higher rank JSJ decomposition we associate with it a 2-orbifold $O_Q$. 
The 2-orbifold $O_Q$ is obtained from the 2-orbifold that is
associated canonically with the 2-orbifold $\hat O_Q$ that has the fundamental group $Q$. If $\hat O_Q$ is an orientable 2-orbifold, we set $O_Q$ to be the
orbifold that is obtained from $\hat O_Q$ by attaching (disjoint) disks to its boundary components. If $\hat O_Q$ is non-orientable, and has no reflection lines
that intersect boundary components,
we set $O_Q$ to be the orbifold that is obtained from $\hat O_Q$ by attaching (disjoint) disks to its boundary components. 
If $\hat O_Q$ has reflection lines that intersect 
boundary components, we set $O^1_Q$ to be the orbifold that is obtained from $\hat O_Q$ by taking out regular neighborhoods of all the unions of
a reflection line that
intersects boundary components together with the boundary components that it intersects. We set $O_Q$ to be the 2-orbifold that is obtained from 
$O^1_Q$ by attaching disks to its boundary components. 

Given a virtually abelian vertex group $VA$ in the higher rank JSJ decomposition, we associate with it a free abelian group $A$ which is the quotient of
$VA$ by the normal subgroup that is generated by the edge groups that are connected to $VA$.

We set the objects in the groupoid $GD3$ to be the sets of 2-orbifolds and free abelian groups that are associated with each higher rank JSJ decomposition in the
higher rank MR diagram, where homeomorphic 2-orbifolds and isomorphic free abelian groups that are associated with the same factor of the higher rank JSJ
are numbered. There are finitely many $m$-collections in the higher rank diagram, so there are finitely many objects in $GD3$.  

In constructing the groupoid $GD1$ we associated arrows with each outer automorphism $\tau in Out(G)$ (theorem 2.1), finitely many arrows with each
outer automorphism. We start associating arrows in the groupoid $GD3$ with 
automorphisms in $Out(G)$ in a similar way. Later we close the set of arrows that are associated with  outer automorphisms, to get a groupoid.

Let $\{\varphi_s\}$
be a weak test sequence of a maximal complexity $m$-collection in the higher rank MR diagram. We say that $\{\varphi_s\}$ is a $durable$ weak  test sequence if
for every $\sigma \in Aut(G)$,  every pair of convergent subsequences of the pair, $\{\varphi_s\}$ and $\{\varphi_s \circ {\sigma}\}$, satisfies conclusion
(5) in lemma 1.13.

Let $\{\varphi_s\}$ be a durable test sequence that factors through a maximal complexity $m$-collection $Coll$, and let $\tau \in Aut(G)$ be an automorphism. 
Suppose that  $\varphi_s \circ {\tau}^{-1}\}$ is a subsequence that factors through a maximal complexity $m$-collection $Coll_1$. Clearly it is a durable test
sequence of $Coll_1$. 

Hence, with every convergent pair
of durable subsequences,  $\{\varphi_s\}$ and $\{\varphi_s \circ {\tau}^{-1}\}$, there are associated automorphisms between 
the QH and virtually abelian vertex groups
in $Coll$ onto the QH and virtually abelian vertex groups in $Coll_1$. We denoted these automorphisms in constructing the groupoid $GD1$ (theorem 2.1) 
$\mu_{\tau}$.

Let $Q$ be a QH vertex group in the higher rank JSJ of the $m$-collection $Coll$, and suppose that it is mapped isomorphically onto the QH vertex group $Q_1$ in
the $m$-collection $Coll_1$ by the map $\mu_{\tau}$. The isomorphism corresponds to a homeomorphism of the orbifold $\hat O_{Q}$ 
onto the homeomorphism $\hat O_{Q_1}$, that maps boundary
components to boundary components. 
With the homeomorphism from $\hat Q_{Q}$ onto $\hat O_{Q_1}$ we can associate a homeomorphism
between the orbifolds $O_Q$ and $O_{Q_1}$. that are obtained by gluing discs to some boundary components or s.c.c.\  in the two orbifolds. 

In a similar way an isomorphism from a virtually abelian vertex group $VA$ in the higher rank JSJ of the $m$-collection $Coll$ 
onto a virtually abelian vertex group $VA_1$ in the higher rank
JSJ that is associated with the $m$-collection $Coll_1$, restricts to an isomorphism between the two associated quotient free abelian groups $A$ of $VA$  and 
$A_1$ of $VA_1$. 

Hence, with every durable weak test sequence $\{\varphi_s\}$, and every (outer) automorphism $\tau \in Aut(G)$, we have associated at least one arrow 
in the groupoid $GD3$. By construction, with each arrow that is associated with an automorphism $\tau$ there is an inverse arrow that is associated with $\tau^{-1}$.
However, the construction does not guarantee that it is always possible to compose arrows. Hence, if $\tau_1 \in Aut(G)$ is 
associated with an arrow from $A$ to $B$,
and $\tau_2$ with an arrow from $B$ to $C$, we associate with $\tau_2 \circ \tau_1$ an arrow from $A$ to $C$, which is obtained as a composition of the
arrows (isomorphisms) from $A$ to $C$ (the objects $A$, $B$, $C$ need not be distinct).

The construction of $GD3$ guarantees that it is a groupoid. However, to get a representation of a finite index subgroup of $Out(G)$ from it, we still need
to prove that the association of arrows between two fixed objects in $GD3$ to outer automorphisms is unique.

\proclaim{Lemma 2.7} Let $A$ and $B$ be two (not necessarily distinct) objects in $GD3$. Let $\tau \in Aut(G)$, and suppose that
there is an arrow from $A$ to $B$ that is associated with $\tau$. Then there is a unique such arrow. Furthermore, up to a composition with an inner automorphism,
the arrow depends only on the class of $\tau$ in $Out(G)$.
\endproclaim

\nfp Let $\tau \in Aut(G)$ be an automorphism with an arrow with a map $\mu_{\tau}$ between the objects $A$ and $B$ in the groupoid $GD3$.
Let $Q$ be a QH vertex group in the higher rank JSJ decomposition that is associated with the object $A$. 
Suppose that $\mu_{\tau}(Q)=Q_1$ is a QH vertex group in the higher rank JSJ decomposition that is associated with the object $B$.

Let $\hat O_Q$ and $\hat O_{Q_1}$ be the 2-orbifolds that are associated with $Q$ and $Q_1$, and let $O_Q$ and $O_{Q_1}$ be the closed 2-orbifolds that
are obtained from $\hat O_Q$ and $\hat O_{Q_1}$ by attaching disks to their boundaries or some s.c.c.\ By killing the generators of all the vertex groups in the
higher rank JSJ decompositions that are associated with the objects $A$ and $B$ that are not
$Q$ or $Q_1$, we get quotients maps $p_Q:H \to \pi_1(O_Q)$ and $p_{Q_1}:H \to \pi_1(O_{Q_1})$, where $\pi_1(O_Q)$ is the (closed) orbifold
fundamental group. Furthermore, the map $\mu_{\tau}$ restricts to a map  $\tilde \mu_{\tau}:\pi_1(O_Q) \to \pi_1(O_{Q_1})$ so that:
$\tilde \mu_{\tau} \circ p_Q =p_{Q_1} \circ \tau$.

The quotient maps $p_Q$ and $p_{Q_1}$ depend only on the higher rank JSJs that are associated with the objects $A$ and $B$, and not on the arrow with the map
$\mu_{\tau}$. Hence, the restriction of the map $\mu_{\tau}$, 
$\tilde \mu_{\tau}:\pi_1(O_Q) \to \pi_1(O_{Q_1})$, depends only on the automorphism $\tau$. Furthermore, the class of $\mu_{\tau}$, i.e., $\tilde \mu_{\tau}$ up to
composition with an inner automorphisms of the (closed) 2-orbifold group, depend only on the class of $\tau$ in $Out(G)$. 
  
Exactly the same argument that we used for orientable 2-orbifolds applies to free abelian  (quotient) groups
that are associated with virtually abelian vertex groups in the
higher rank JSJ decompositions. 

\line{\hss$\qed$}

Given the groupoid $GD3$, we look at the fundamental group w.r.t.\  one of the objects in $GD3$. The fundamental group is some finite index subgroup
of $Out(G)$, and the groupoid gives a homomorphism of the fundamental group to the direct sum of the mapping class groups of the (closed) 2-orbifolds that
are associated with the QH vertex group in the higher rank JSJ that is associated with the object, and with the direct sum of the general linear groups
which are the outer automorphism groups of the free abelian (quotient) groups that are associated with the virtually abelian groups in this higher rank
JSJ. This completes the proof of theorem 2.6.

\line{\hss$\qed$}

Theorem 2.6 shows that with the groupoid $GD1$ it is possible to associate a homomorphism from a finite index subgroup of $Out(G)$ into the direct sum
of finitely many mapping class groups of closed orbifolds and finitely many general linear groups. Hence, given the groupoid $GD1$, over each value of
this homomorphism (i.e., the image of an element from $Out(G)$) there is a finite fiber that is encoded by $GD1$.

In the case of a hyperbolic group, the kernel of the homomorphism from a finite index subgroup of $Out(G)$ onto the direct sum of the mapping class
groups of the 2-orbifolds in the JSJ decomposition of $G$ is f.g.\ virtually abelian (see [Se7] and [Le]). In the case of hyperbolic groups,
the analysis of this kernel is equivalent to the
analysis of $Out(G)$ in case the JSJ of the group has no QH vertex groups.

In the case of HHG we constructed the groupoid $GD2$ in case the higher rank JSJ has no QH and no virtually abelian vertex groups. We were not able to prove that
in this case $Out(G)$ is f.g.\ and virtually abelian, but we do prove that it is locally virtually nilpotent.  

\proclaim{Theorem 2.8} With the assumptions of theorem 2.5, suppose that the higher rank JSJ decomposition of $G$, contains no QH and no virtually abelian vertex groups.
Then $Out(G)$ is locally virtually nilpotent (with uniformly bounded nilpotency class).
\endproclaim

\nfp Let
 $M=<\tau_1,\ldots,\tau_{\ell}>$ be a f.g.\ subgroup in $Out(G)$. By theorems 2.4 and 2.5 the number of elements, $\nu^j_{\tau} \in Z^n$, that are associated 
with each of the elements
$\tau_i$, $i=1,\ldots,\ell$, is finite. Let $b$ be the maximal norm of these elements.

First, suppose that all the cover resolutions in  the $m$-collections in the higher rank MR diagram of $G$ have at most two levels (the assumption in theorems 2.1
and  2.4). 
The number of elements, $\nu^j_{\tau} \in Z^n$, that are associated with all the elements in a ball of radius $r$ in $M$
(w.r.t.\ its fixed set of generators) is bounded by $(rcb)^n$ for some universal constant $c$ (that depends on the norm). By the discreteness of the action of $G$ on the
product $X$, it follows that the number of elements in a ball of radius $r$ in $M$ is at most $c_1(rcb)^n$, hence, it is at most polynomial. By Gromov's polynomial growth
theorem it follows that $M$ is virtually nilpotent (with a universal bound on its nilpotency class).

Now, suppose that the higher rank MR diagram of $G$ is a general single ended diagram. If the edge groups in the JSJ like decompositions can be conjugated 
to the terminal level in the higher rank MR diagram,  the argument that was used in case there are at most two levels in each resolution in the higher rank diagram is
valid. Otherwise, for each edge group that can not be conjugated into the bottom level, there can be at most  $rc_2b$ elements that belong to the edge group and can serve as
conjugating elements for all the automorphisms of a ball of radius $r$ in $M$. When the conjugating elements that are associated with the edge groups are fixed, the vertex
groups in the higher rank JSJ can have at $c_3$ values. Hence, the number of elements in a ball of radius $r$ in $M$ is bounded by
$c_3c_4(rc_2b)^n$ for some constant $c_4$ that depends on the norm. By Gromov's polynomial growth theorem, $M$ is virtually nilpotent (with a universally bounded
nilpotency class).

\line{\hss$\qed$}

Theorems 2.6 and 2.8 naturally raise the following questions, that we were unable to answer:

\proclaim{Question 1} Can the number of morphisms that are associated with each automorphism in $Out(G)$ in theorems 2.1, 2.4 and 2.5 be uniformly bounded?
we only proved that it is finite. Note that if there is only a single morphism that is associated with each outer automorphism in $Out(G)$ and a source object in the
groupoid, then in the case of the groupoid $GD1$, it is possible to obtain a homomorphism from $Out(G)$ into direct products of mapping class groups and
outer automorphism groups of virtually f.g.\ abelian groups, which is a direct generalization of the structure of $Out(G)$ in the case of  hyperbolic groups 
and of theorems 1.9 and 1.14.
\endproclaim

\proclaim{Question 2} In theorem 2.8 we proved that if the higher rank JSJ contains no QH and no virtually abelian vertex groups, then $Out(G)$ is
locally virtually nilpotent with a universally bounded nilpotency class. Is $Out(G)$ f.g.? Is it virtually abelian?
\endproclaim

\vglue 1.5pc
\centerline{\bf{\S3. A higher  rank JSJ decomposition of (some) HHG }}
\medskip

In the previous section we used the higher rank MR diagram that was constructed in [Se8], to construct the higher rank JSJ decomposition of a group $G$
that acts discretely and cocompactly on a product of (finitely many) hyperbolic spaces. With the higher rank JSJ and $Out(G)$
we associated two groupoids $GD1$ and $GD2$ that contain information on the algebraic structure and the dynamics of individual automorphisms in $Out(G)$.
We further constructed a homomorphism from a finite index subgroup of $Out(G)$ into the direct sum of mapping class groups of
orbifolds that are associated with QH vertex groups in the higher rank JSJ,  and general linear groups that are associated with the virtually abelian vertex groups
in the higher rank JSJ decomposition (theorem 2.6). 

In [Se8] we first constructed a higher rank MR diagram for the automorphism group of a group that acts discretely and cocompactly on a product of hyperbolic
spaces (theorem 2.8 in [Se8]), and then modified and generalized the construction of the higher rank MR diagram to colorable HHGs that satisfy a weak acylindricity
condition. 
These conditions do not hold for all HHG (e.g.\ Burger-Mozes groups), but they both hold for the mapping class groups, due to seminal works of
Bowditch [Bo] and Bestvina-Bromberg-Fujiwara [BBF1].
    
In section 3 of [Se8], we analyze resolutions of automorphisms of an HHG $G$, and constructed a higher rank MR diagram,  under two possible assumptions.  
Under both assumptions it is possible to
generalize the analysis of automorphisms of an HHG that acts on a product that was applied in the first two sections. We start with the stronger assumption
from section 3 in [Se8].

\proclaim{Theorem 3.1} Let $G$ be a colorable HHG. 
Let $H<G$ be characteristic finite index subgroup, 
such that domains in an orbit of $H$ are pairwise transverse,
and the (finitely many) actions of $H$ on the quasi-trees of metric spaces that are constructed from the action of $H$ on these orbits
via the [BBF1] construction, for all the quasi-trees of metric spaces that are not quasi-isometric to a real line,   are weakly acylindrical.

By theorem 3.3 in [Se8], with these assumptions we can associate with $Aut(G)$ a higher rank MR diagram. Suppose that all the virtually
abelian decompositions along the cover resolutions in the $m$-collections in this diagram do not contain edges with finite (nor trivial) edge groups.

Then it is possible to associate with $Aut(G)$ a higher rank JSJ decomposition, similar to the one that was constructed in case $G$ acts on a product of hyperbolic
spaces in the previous section (theorem 2.5). With the higher rank JSJ it is possible to associate two groupoids, 
similar to either the groupoid $GD1$ that was constructed in theorem 2.1,
or $GD2$ that was constructed in theorem 2.4, precisely as we associated these groupoids with HHGs that act on products in theorem 2.5.

Finally with the higher rank JSJ decomposition we can associate a homomorphism from a finite index subgroup of $Out(G)$ into the direct sum of mapping class 
groups of finitely many orbifolds direct sum with finitely many general linear groups, precisely as it was constructed in theorem 2.6. If the higher rank JSJ decomposition
contains no QH and no virtually abelian vertex groups in its JSJ like decompositions, then $Out(G)$ is locally virtually nilpotent  (with uniformly bounded nilpotency class).
\endproclaim

\nfp 
Under the assumptions of the theorem, it is possible to associate with $Aut(G)$ a higher rank MR diagram (theorem 3.3 in [Se8]).
Having a higher rank MR diagram, and the weakly acylindrical action of $H$ on the (finitely many) quasi-trees of metric spaces, that are constructed from
the action of $H$ on its orbits using [BBF1],  
all the analysis of automorphisms of HHGs that act on product spaces, as it appears in the first two
sections of this paper, generalizes to HHGs that satisfy the assumption of the theorem. In particular, the conclusions of theorems 2.1, 2.4, 2.5, 2.6 and 2.8  
remain valid,
and we can associate with $Aut(G)$ a higher rank JSJ decomposition, the groupoids $GD1$ and $GD2$, and the homomorphism that is constructed in theorem
2.6.

\line{\hss$\qed$}

The second possible assumption in section 3 of [Se8] is weaker and requires finer analysis and somewhat weaker notion of a cover resolution, that is
called $hybrid$ $resolution$ in theorem 3.4 in [Se8], in order to
construct a higher rank MR diagram. 

Recall that a hybrid resolution consists of a pair of resolutions.    
The first  resolution in each pair, is a resolution of the
ambient group $H$, and it is constructed from actions of $H$ twisted by automorphisms from $Aut(G)$, 
on the associated projection complexes $P_K^j$, $1 \leq j \leq m$ (see [BBF1] for the construction of projection complexes). 

The second part in each pair is a resolution of the terminal limit group of the first resolution. The second resolution
 is composed from finitely many resolutions of
subgroups of the terminal group of the first resolution, which are the intersections of the terminal group with (set) stabilizers of domains in $X$.

By theorem 3.4 in [Se8], the higher rank MR diagram consists of finitely many $m$-collections of hybrid resolutions.
Every automorphism in $Aut(G)$ factors through at least one of the finitely many $m$-collections of hybrid resolutions. Note that an automorphism factors through  
an $m$-collection of hybrid resolutions if the action of $H$ on the HHS space $X$ factors through all the $m$ hybrid resolutions in the $m$-collection according
to definition 2.6 in [Se8].

The next theorem uses this higher rank MR diagram of hybrid resolutions, that was constructed in theorem 3.4 in [Se8],
and  applies the analysis of automorphisms that appears in the first two sections of this paper, to generalize the construction of the higher rank
JSJ decomposition and the groupoids and homomorphism that are associated with it.

Note that both colorability and the weak acylindricity assumptions in theorem 3.2, that are the assumptions in theorem 3.4 in [Se8],  
are satisfied by the mapping class groups of surfaces.

\proclaim{Theorem 3.2} Let $G$ be a colorable HHG, and let $H<G$ be a characteristic finite index subgroup,
such that in each orbit under the action of $H$ domains are transverse.
We further assume that $H$ acts weakly acylindrically on the (finitely many) projection complexes that are constructed from the actions of $H$ on these orbits
of domains (using [BBF1]),
and the set stabilizers of each domain acts weakly acylindrically on the domain. 

In theorem 3.4 in [Se8], with these assumptions we can associate with $Aut(G)$ a higher rank MR diagram. 
Let $m$ be the number of orbits of projection spaces under the action of $H$. Each $m$-collection in the
higher rank MR diagram, contains $m$ $hybrid$ cover resolutions. Each hybrid cover resolution consists of a resolution of an ambient limit
quotient of $H$, followed by finitely many resolutions of quotients of f.g.\ subgroups.

Suppose that all the virtually
abelian decompositions along the hybrid resolutions in this diagram do not contain edges with finite (nor trivial) edge groups.

Then it is possible to generalize the conclusions of theorems 2.1, 2.4, 2.5, 2.6 and 2.8 to HHG that satisfy these assumptions.
It is possible to associate with $Aut(G)$ a higher rank $hybrid$ JSJ decomposition, groupoids similar to $GD1$ and $GD2$ from theorems 2.1 and 2.4, 
and a homomorphism
similar to the one that was constructed in theorem 2.6. Furthermore, if the hybrid JSJ like decompositions in the hybrid higher rank JSJ decomposition contains no QH and no
virtually abelian vertex groups, $Out(G)$ is locally virtually nilpotent with uniformly bounded nilpotency class.

\endproclaim

Note that like in the construction of the higher rank MR diagram in theorem 3.4 in [Se8], with a general HHG we associate a $hybrid$ JSJ
decomposition, which generalizes and is somewhat different than the higher rank JSJ in case of a group that acts on a product
(theorems 2.1, 2.5 and 2.6). The precise definition and construction of the hybrid JSJ for a colorable HHG  that satisfies our weak acylindricity assumptions appears in the
sequel. 

\nfp Let $m$ be the number of orbits of domains under the action of $H$. By theorem 3.4 in [Se8] under the assumptions of the theorem, it is possible to associate
with the action of $H$ on the HHG $X$ and its domains, a higher rank MR diagram that contains finitely many $m$-collections of cover hybrid resolutions 
(see theorem 3.4 in [Se8] for these notions and statement).

Each (cover) hybrid resolution in an $m$-collection from the higher rank MR diagram 
has a top and a bottom part. The top part is a (single) cover resolution that is obtained
from the action of $H$ on a projection complex that is constructed from the corresponding orbit of domains under the action of $H$. The limit group that is
associated with such a resolution is a (cover of a) limit quotient of the ambient group $H$.

Each bottom part contains finitely many resolutions of some  (covers of) f.g.\ subgroups of the  terminal limit group of the top resolution. 
These resolutions are constructed from the action of the terminal limit group of the top resolution on the quasi-trees of metric spaces that are constructed
from the action of $H$ on the orbits of the domains using [BBF1]. Each resolution from the finite set that is associated with the bottom part
of a hybrid resolution is constructed from the action of the (set) stabilizer of a 
domain on the domain. 

We set the complexity of a hybrid resolution to be the complexity of its top resolution 
(the one that is obtained from the action of $H$ on the corresponding projection complex) followed by the tuple of complexities of the resolutions in its
bottom part, where these last complexities are ordered in a non-increasing order. 

On the set of complexities of hybrid resolutions we can define an order which is basically a lexichographical order. Given the complexities
of two hybrid resolutions, we first compare the complexities of their top resolutions, and if these are equal we compare their tuples of complexities of
their bottom resolutions, where the comparison of the two tuples is lexichographical.

Given an $m$-collection of hybrid resolutions, we define its complexity to be the $m$ complexities of its $m$ hybrid resolutions ordered in a non-increasing order.
On the set of complexities of $m$-collections of hybrid resolutions we naturally define the lexichographical order.

In theorem 3.4 in [Se8] we constructed a higher rank MR diagram with $m$-collections of hybrid resolutions to study the structure of the
automorphism group of a colorable HHG that satisfies the weakly acylindrical assumption.
Once we defined the complexity of an $m$-collection of hybrid resolutions, the same arguments that were used to prove theorems 2.1,2.4,2.5, 2.6 and 2.8  in the
case of a group that acts on a product of  hyperbolic spaces, and $m$-collections of ordinary resolutions, can be used to prove the analogous
 statements for colorable HHG that satisfy the weak
acylindricity condition. 

However there is a major difference between the higher rank JSJ decomposition that we constructed in the case of a group that acts on a product of
finitely many hyperbolic spaces, 
and the higher rank JSJ  of a general colorable HHG that we call a $hybrid$ higher rank JSJ. In the case of a group that acts on a product, the higher rank
JSJ is an $m$-collection of  virtually abelian decompositions of some quotients of a finite index subgroup $H$ of the ambient group $G$ that acts on a
product of hyperbolic spaces.

In the case of a colorable HHG $G$ that satisfies the weak acylindricity condition, the $hybrid$ higher rank JSJ contains $m$ decompositions of quotients
of a finite index subgroup $H$ of $G$ or alternatively f.g.\ groups that contain these quotients, see the proof of theorem 2.5 for these two possible options), 
but these decompositions need not be virtually abelian in general. Indeed, the part in each decomposition that is
constructed from the top resolution, i.e., from the action of $H$ on the projection complex that is associated with a given orbit, is a virtually abelian
decomposition.  

However, the parts in each decomposition in the hybrid higher rank JSJ  that are constructed from the action of (set) stabilizers of domains on their domains, 
are in general not virtually abelian. These parts are virtually abelian modulo the pointwise quasi-stabilizers of these domains. Hence, edge stabilizers and 
virtually abelian and QH vertex stabilizers in these parts of the hybrid higher rank JSJ contain (covers of) the (pointwise) stabilizers of the 
corresponding domains.
Therefore, these parts of decompositions in the the hybrid JSJ are not virtually abelian decompositions.

\line{\hss$\qed$}


\vskip 2em

\end